\documentclass[12pt]{article}
\usepackage{a4,verbatim,amsmath,amssymb}
\usepackage[isolatin]{inputenc}
\usepackage{epsfig}
\usepackage[all]{xy}

\input math1.def
\usepackage{amsthm}

\newtheorem{theor}{Theorem}[section]

\newtheorem{lemma}[theor]{Lemma}
\newtheorem{cor}[theor]{Corollary}
\newtheorem{prop}[theor]{Proposition}

\def\proof{\goodbreak\noindent{\sc Proof. }\nobreak}

\def\endproof{\par\nobreak\hbox to \hsize{\hfil\vrule width 5pt height
5pt}\goodbreak\vskip 3pt}

\allowdisplaybreaks \numberwithin{equation}{section}

\begin{document}
\title{Green's  Mapping and Julia Sets}
\author{Ilia Binder\thanks{Partially supported by an NSERC Discovery grant. }, \,\,  Paul F. X. M\"uller\thanks{Partially supported
    by the Austrian Science foundation (FWF) Pr.Nr. FWFP28352-N32.}  \,\,and Peter Yuditskii\thanks{Partially supported
    by the Austrian Science foundation (FWF) Pr.Nr. FWFP28352-N32.}}
\date{\today}
\maketitle

\section{Introduction}

Fatou's theorem published in 1906 asserts that  any bounded analytic function $f : \bD \to \bC $ admits  almost sure radial limits.  Specifically,
there exists $ M_f \subset \bT$ such that $m(M_f) = 0 ,$ and  
\begin{equation}\label{intro-F}
  \lim_{r \to 1} f(r\zeta) \quad\text{exists for }\quad \zeta \in \bT \setminus M_f .
  \end{equation}
  Thus Fatou's theorem
led to the search for improved boundary behavior of analytic functions $ f : \bD \to \bC $.
In this direction  A. Beurling in 1933, 
proved that the Dirichlet condition
 $$ \int_{\bD} |f'(z)|^2 dA(z) < \infty $$
 implies the existence of   an exceptional set $M_f$,  with vanishing logarithmic capacity,  satisfying  
\begin{equation}\label{RV}
\int_0^1  |f'(r\zeta)| dr  < \infty \quad\text{ for }\quad \zeta \in \bT \setminus M_f .
\end{equation}
In 1993, J. Bourgain   showed that for every bounded  analytic function $f$
on the unit disk and any interval  $ I \subset \bT $ there exists a set of
directions $ B_f \subset I$ satisfying 
\begin{equation}\label{B}
  Hdim (B_f) = 1 \quad \text{and} \quad \int_0^1  |f'(r\zeta)| dr  < \infty , \text{ for }  \zeta  \in B_f  .
\end{equation}
In this paper we investigate a class of analytic functions (the Green's mappings) on the unit disk for which
the following
quite restrictive boundary convergence condition holds true:
\begin{equation}\label{RV}
  \int_0^1  |f^{\prime \prime}(r\zeta)| dr  < \infty ,  \quad\text{ for at least one }
  \quad \zeta \in \bT.
\end{equation}
For background and motivation we first review  condition \eqref{RV} for 
conformal functions and universal covering maps, which are closely related to Green's mappings. 
\subsection{Conformal Functions  and universal covering Maps}
We begin with conformal functions on the unit disk. Let
$\varphi : \bD \to \bC $  be  analytic and injective. 
In 1971 J.M. Anderson~\cite{anderson}  conjectured that  
 there exists $\zeta \in \bT$ such that
\begin{equation}\label{conformal}
   \int_0^1 |\varphi^{''} (r\zeta) |dr < \infty.
 \end{equation}
In 1997 the  proof of 
  Anderson's conjecture  was obtained in  \cite{jones-mueller-1997}.
  Conceptually,  Anderson's problem 
 is treated as a 
 problem on the radial variation of the real harmonic Bloch function $ b = \log |\varphi '| $,
 satisfying 
\begin{equation}\label{bloch}
   \sup_{z\in \bD} |\nabla b(z)|(1-|z|) < \infty , 
   \end{equation}
 It was shown that  \eqref{bloch} implies that there exist  $\zeta \in \bT$ such that
\begin{equation}\label{JM}
   \int_0^1 |\nabla b(r\zeta)| e^{b(r\zeta)} dr < \infty,
 \end{equation}
Since $2|\varphi ^{''}(z)| =  |\nabla b(z)| e^{b(z)}$, \eqref{JM} gives \eqref{conformal}, in which case we say
that $ \zeta $ is a good direction for $\varphi $.
The proof in \cite{jones-mueller-1997} shows that the Hausdorff dimension of the good directions for  $\varphi $ equals $1 . $

We continue  our introductory  discussion by reviewing variational estimates for
universal covering maps onto planar domains:
Fix a closed set $ E \subset \bC$.
By the uniformization theorem, if the cardinality of $E$
is $\geq 2$, the universal covering surface  of $ \bC \setminus E $ is conformally equivalent
to $\bD . $  We let
$$P_E : \bD \to  \bC \setminus E $$
denote the corresponding  universal
covering map. For two important special classes of closed sets  $  E $, 
variational estimates for  $P_E $ are readily  obtained  by reduction
to  \cite{jones-mueller-1997}:  
\begin{enumerate}
\item 
If $E$ is connected then  $\bC \setminus E$ is simply connected and
$P_E : \bD \to  \bC \setminus E $ is bijective. Hence, by \cite{jones-mueller-1997}, $P_E$ satisfies the estimates of
Anderson's conjecture \eqref{conformal}.
\item
If more generally  $E$ is a uniformly perfect set then,  
by a theorem of Ch. Pommerenke~\cite{pommerenke-unifperf},  
$b_E = \log |P_E '|$ is a harmonic Bloch function on $ \bD . $ Hence again,
by \cite{jones-mueller-1997}, we have \eqref{JM}, and consequently \eqref{conformal} holds true, with $\varphi$ replaced by $P_E$.  If  $E$ is  uniformly perfect, the set of good directions 
is of Hausdorff dimension $=1. $ 
\end{enumerate}

In  \cite{jones-mueller-1999} it was shown 
that variational estimates, analogous to \eqref{conformal},
holds true for  universal covering maps  in full generality:   
Despite the fact that  for {\em general}  universal covering maps,  $\log |P_E^{'}| $ need not be a Bloch function, we showed in   \cite{jones-mueller-1999} that there exist countably many
$\psi\in [0,2]$ such that

\begin{equation}\label{ucova-int1}
  \int _{0}^{1}|P_E''(re^{i\pi\psi})|dr<\infty.
\end{equation}
(In case where   $ E = \{ 0,1 \} ,$   the set of good directions for $P_E$ is at most countable, showing that the above assertion is  best possible.)
%
 
%

\subsection{Green's mapping} 
The notion of the Green's mapping associated to a complex domain
is closely related to that of  universal covering maps.
We refer here to 
the classical works
Arsove and Johnson \cite{arsove-johnson-1970},
and  Pommerenke~\cite{pommerenke-kap, pommerenke-1977}.

We fix a non polar domain $ \Omega\subset \mathbb{C}$.  
Let $g:\Omega\times\Omega\rightarrow (-\infty,\infty]$
denote the  Green's {\em function} for $\Omega .$
Let $z_0$ be a fixed point in $\Omega$.
The level lines of the Green's function $g(\cdot,z_0)$ will be called Green's lines. Arcs in $\Omega$ not passing through singular points and such that the tangent of the arc and the gradient of the Green's function are parallel in each point are called Green's arcs.
Let  $\tilde{\Omega}$ be the simply connected subdomain of $\Omega$ consisting  of all Green's arcs with starting point $z_0$, (the so called central Green's arcs.)
We let  $H \sb \mathbb{D}$ demote the   {\em Green's fundamental domain} of $\Omega$, which, we remind the reader, is a radial slit domain,  called sometimes  hedgehog domain.  
 
We let
$$T: H \to \tilde{\Omega}$$
denote   {\em Green's mapping}  of $\Omega . $  
It is well known that the  {\em Green's fundamental domain}  $H$, reflects the geometry of $\Omega$ and its boundary.
For instance for Denjoy domains, i.e.,
$\partial \Omega \subset \bR$, 
V. Andrievski~\cite{andrievski} proved that $H$ is a John domain if and only if $\partial \Omega$
is uniformly perfect.

In March 1999, the first named author (Binder) posed the problem of showing that  
a ``good direction'' $\psi\in [0,2]$ exists, for any  Green's mapping $T:H\rightarrow\tilde \Omega$, i.e., 
\begin{equation}\label{binder}
  \int\limits_0\limits^{1} |T''(re^{i\pi\psi})|dr <\infty, \quad\text{ for at least one }
  \quad \psi\in [0,2], 
\end{equation}
and delineated a path connecting  \eqref{binder} to \eqref{ucova-int1}.
Presently this problem is open  even in the special case where   $\partial \Omega $ is a uniformly perfect subset of the real line. In this paper we obtain a positive solution
when $\Omega = \overline{\bC} \sm E_0$  where $E_0 \sb \bR $ is  the Julia set of an expanding quadratic polynomial.

\subsection{Main Result}

A primary example of a Denjoy domain $\Omega = \overline{\bC} \sm E_0$  with a uniformly perfect boundary $E_0\subset\bR$ of {\em vanishing  Lebesgue measure} is the Julia set of an expanding quadratic polynomial $P(z)=z^2-\lambda$ with $\lambda>2$. The Julia set $E_0$ is defined as a complement to the basin of attraction of infinity
$$
E_0=\bC\setminus\{z\in\bC: P^{\circ n}(z)\to\infty\}
$$
where   $P^{\circ n} $ is the $n-$fold iterate of $P$, $P^{\circ (n+1)}=P\circ P^{\circ n}$.

Despite its elementary nature, the bifurcation diagram of the family $P(z)=z^2-\lambda$ with $\lambda\in\bC$,  has an extremely rich and complex structure. It was Mandelbrot's computer-generated image of this diagram that aroused great interest in problematic's in general \cite[Chap. 2, Sect. 7]{eremenko-lyubich}
. Also, Jacobi matrices generated by iterations of quadratic polynomials are the leading model in inverse spectral theory.  Starting from the pioneering work of Bellisard--Bessis--Moussa \cite{BBM82}, they have attracted much attention, providing almost periodic operators with singular continuous spectrum; cf. \cite{BGH82, BGH85, PVY06}, see \cite{ELYps} for a recent development. In \cite[Problem 1.4]{VanAssche}, Simon asked about the fine structure of zeros of orthogonal polynomials for singular measures,  for a partial answer see \cite{ELYrl}. The same problem deals with asymptotics for Chebyshev polynomials (monic polynomials with the least deviation from zero on compact subsets of the real axis
\cite{sodin-yuditski-92}). In both cases  polynomials related to Julia sets provide the best understood  examples.

Thus it is only natural to test the first named author's (Binder's) conjecture  for real Julia sets. To state our result we need some extra notations. 

\begin{figure}   
    \begin{center}
\includegraphics[scale=0.5]{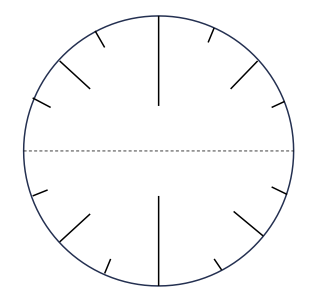}
\end{center}
\caption{Green's Fundamental domain $H$}
  \label{fig:2}
		\end{figure}   

                Green's fundamental domain $H$ associated to a Denjoy domain is symmetric with respect to the real axis,  see Figure \ref{fig:2}. If moreover the boundary of a  Denjoy domain is
                a Julia set, conjecture  \eqref{binder}  can be effectively analyzed due to special arithmetical features of $H$. For a given $\lambda>2$ there exists $ a=a(\lambda) > 0 $ such that  
the  \textit{symmetric} domain $H=H(a)$ is obtained by removing from the lower semi-disc $\bD_-= \{ \z: |\z|<1, \Im\z<0\} , $
the following sequence of  radial intervals 
$$\{ \zeta=e^{-\pi(h+i\psi_{km})}:  \psi_{k,m} =  k /  2^m,\  h\in (0, a /  2^m]\},   $$
where $m\in \bN $,     $k$ is an odd integer, $k \le 2^m$. 

It is convenient to use dyadic expansion for an arbitrary point in the lower semi-circle 
$$
\bT_-=\{\zeta=e^{-i\pi\psi},\quad   \psi = \sum_{k = 1}^\infty \varepsilon_k 2^{-k},
  \quad\text{where}\   \varepsilon_k \in \{0, +1\}\}. 
  $$ 
  We say that $ \psi  \in F_N $ if $\varepsilon = ( \varepsilon_k) $ satisfies the following conditions:
  For every $ k \in \bN $ the hypothesis
  $\varepsilon_k =\varepsilon_{k+1} =   \dots = \varepsilon_{k+N}  $ implies that
  $ \varepsilon_{k+N+1} \not = \varepsilon_{k+N}  $.

  \begin{theor}[Main Theorem]\label{main} 
  Let $P(z)=z^2-\lambda$ be an expanding quadratic polynomial with a  real Julia set $E_0$.  Let $H$ be the 
  Green's fundamental domain associated to the domain  $\Omega=\overline{\bC}\setminus E_0$, $T:H\to\tilde \Omega$. Let $\lambda>2+\sqrt{2}$, then
   for $ N \in \bN$, there exists $ C_N < \infty $ such that  for any  $\psi \in F_N $,
  $$ \int _{0}^1 |T^{''} (e^{ -i\pi  \psi}r) | \cdot d r \le C_N $$
  and
$${\rm Hdim}(F_N) > 1 - c/N , $$ 
for some  $ c> 0 $ (independent of $N$).  
   \end{theor}

   Thus the first named author's conjecture
   is confirmed in the stated case. Moreover the collection of ``good directions''  $F_\sigma=\cup_{N\in\bN} F_N$ has 
Hausdorff dimension $= 1$.

\subsection{Structure of the paper}
We exploit comb domain techniques \cite{eremenko-yuditski-2012} in  combination with ideas  originating with  Carath\' eodory kernel theorem
(similar to the derivation of Loewner's chain equation).

Let $\xi$ be the positive fixed point of $P$,  $P(\xi)=\xi$. The interval $[-\xi,\xi]$ is the smallest interval containing the Julia set $E_0$. Green's mapping $T$ acts from $H$ to
$\overline \bC\setminus [-\xi,\xi]$. Let $G$ be the inverse map $G: \overline \bC\setminus [-\xi,\xi]\to H$. Note that
 $|G(z)| = 1 $ for $z\in E_0$ and   $|G(z)| <  1 $ for $z\notin E_0$. That is,
    $$- \log |G(z) | =g(z)=g(z,\infty)$$
    is just the Green's function of $\overline{\bC} \sm E_0  $ with respect to infinity. In this connection $G(z)$ is called the complex Green function. It  is multivalued (character-automorphic)
    in  $\overline \bC\setminus E_0$.

From the point of view  of Complex Dynamics  $G$  is called B\"ottcher function \cite[p. 34-35]{carleson-gamelin}. In a vicinity of infinity it conjugates the dynamical system generated by $P(z)$ with the  dynamics given by  $w\mapsto w^2$. That is, 
\begin{equation}\label{eq25-1}
G(P(z))=G(z)^2.
\end{equation}

Another special function deals with the fixed point $\xi$. In the vicinity of this point dynamics of $P(z)$ is equivalent to the multiplication by $\rho$ with $\rho=P'(\xi)$. In this case there exists an entire function $F(z)$ such that
\begin{equation}\label{eq25-2}
P(F(z)+\xi)= F(\rho z)+\xi,\quad F(0)=0,\ F'(0)=1.
\end{equation}
This is the so called Poincare equation, we call $F(z)$ the Poincare function.

Both claims \eqref{eq25-1}, \eqref{eq25-2} are classical in the Theory of Dynamical systems and they are valid in  much more general settings, see survey \cite{eremenko-lyubich}
and in particular
\cite{grabner}. 
In the Second Section we provide a self-contained presentation of these identities,  based on comb domain techniques. Thus we define   domains $\Pi_0$, $\Pi_F$ and corresponding conformal mappings
$$
\theta:\bC_+\to\Pi_0,\quad \Psi_F:\bC_+\to\Pi_F,
$$
which are related to the B\"ottcher and Poincare functions  by
 $$   G (z)  =  e^{ i\pi  \theta ( z )} ,  \ z \in  \bC_ +  , \quad F(w ) +\xi =  \xi \cos ( \pi \Psi_F (w) )), \ w\in\bC_+.
 $$

In the last subsection of Section 2 we introduce a superposition of the above mappings
\begin{equation*}\label{}
\Delta:\bC_+\to \Pi\quad \text{such that}\quad e^{i\pi\Delta(z)}= e^{i\pi\theta(F(z))},
\end{equation*}
and give a detailed clarification of its meaning. 
As a consequence of \eqref{eq25-1} and   \eqref{eq25-2}, the comb domain $\Pi$ is self-similar,
$2 \cdot \Pi = \Pi$,  and moreover
\begin{equation}\label{eq25-3}
2\Delta(z)=\Delta(\rho z).
\end{equation}

In the second preparation part (Section 3) we introduce a collection of  Cantor type sets $F_N$, which are invariant in the following sense. Let
$$
E_0^N=\{T(e^{-i\pi\psi}):\ \psi\in F_N\}\subset E_0=\{T(e^{-i\pi\psi}):\ \psi\in [0,1)\}
$$
Then 
\begin{equation}\label{eq25-4}
P: E_0^N\to E_0^N
\end{equation}
This observation is important in what follows.
We estimate from below the  Hausdorff dimension of $F_N$, proving the second statement of the Main Theorem.

In Section 4 we prove our Main Theorem.
Let  $s_{\psi} $ be the vertical path defined by 
$$ s_{\psi} (h) =   ( -\psi +ih), \quad 0\le h < a.  $$
We assume that the path
$  R_{\psi} = 
    e^{ i\pi  s_\psi }\subset H $
    terminates on $\bT$  (as $h=0$), i.e., $\psi\not = k/2^n$.
 In this case the following paths, both  in $\bC_+$, are well defined
 $$ \gamma_\psi (h) =  \theta^ {-1} (s_{\psi} (h)), \quad h < a ,\quad   \Gamma_\psi (h) =  \Delta^ {-1} (s_{\psi} (h)), \quad h < a  .$$  
  We prove the following series of reductions
$$
\begin{matrix}
\int_0^1|T''(r e^{-i\pi\psi})| d r <\infty
\\
\Uparrow\\
\int_{\gamma_\psi}\left|\frac{\theta''(z)}{\theta'(z)^2}\right||dz|<\infty\\
\Uparrow\\
\int_{\Gamma_\psi}\left|\frac{\Delta''(w)}{\Delta'(w)^2}\right||dw|<\infty
\end{matrix}
$$
We put  
$$  \Gamma_{\psi, n} = \{w \in \bC_+ :  \Delta(w) = - \psi +ih : h \in (a/2^{n +1}, a/2^{n }) \} , $$
 and in view of \eqref{eq25-3} we have
$$
\int_{\Gamma_\psi}\left|\frac{\Delta''(w)}{\Delta'(w)^2}\right||dw|=
\sum\frac {2^n}{\rho^n} \int_{\rho^n\cdot \Gamma_{\psi,n}}\left|\frac{\Delta''(w)}{\Delta'(w)^2}\right||dw|.
$$

Let $\eta$ be positive  such that $P(\eta)=-\xi$. Our restriction on $\lambda$, $\lambda>2+\sqrt{2}$, is equivalent to $\eta>1$. Using the invariance property \eqref{eq25-4} and a kind of Harnack's inequality, we show that
$$
\int_{\rho^n\Gamma_{\psi,n}}\left|\frac{\Delta''(w)}{\Delta'(w)^2}\right||dw|\le C\frac{\rho^n}{2^n \eta^n}
$$
  uniformly in $\psi\in F_N$.

Therefore
$$
\int_{\Gamma_{\psi}}\left|\frac{\Delta''(w)}{\Delta'(w)^2}\right||dw|\le C\frac{\eta}{\eta-1}<\infty,
$$
  and Theorem \ref{main} follows.


\section{Preparation I }

 \subsection{The Comb Domain Representation of  the B\"ottcher function for an expanding quadratic polynomial $P(z) = z^2  -\lambda . $}

Let $ a > 0 $.   Given $ u, v \in \bC $ we  denote  $[u,v] = \{ u + t(v-u): 0\le t \le 1 \}$.
The comb domain $\Pi_0(a)$ is obtained by removing from the vertical half strip $ \{ \psi +i \phi : \psi \in (-1, 0 ) , \phi > 0 \} , $
the following sequence of  slit-intervals
$$ [  \psi_{k,m} , \psi_{k,m} +i\phi_{k,m } ] , \quad m\in \bN ,  \quad k \in \bN, \quad k \le 2^m , k \quad \text{odd} ,   $$
where
$$  \psi_{k,m} =  -k /  2^m , \qquad \phi_{k,m } =  a /  2^m . $$
Thus
$$
\Pi_0(a) = \{ \psi +i \phi : \psi \in (-1, 0 ) , \phi > 0 \} \setminus \bigcup_{m \in \bN }  \bigcup _{ k \le 2^m ,\  k \  \text{odd} }   [  \psi_{k,m} , \psi_{k,m} +i\phi_{k,m } ]
$$
  Let $\widetilde{\Pi_0(a)}$ be the reflection of  $\Pi_0(a) $ across the line   $-1 + i(0, \infty)   $ and form the comb domain 
  $$ 2 \cdot\Pi_0(a) = \widetilde{\Pi_0(a)} \cup \{-1 + i(a, \infty) \} \cup \Pi_0(a) . $$
  We remark that $ 2 \cdot\Pi_0(a)$ was  obtained by removing vertical slit intervals from the half strip
    $ \{ \psi +i \phi : \psi \in (-2, 0 ) , \phi > 0 \} . $ Accordingly the base interval of   $ 2 \cdot\Pi_0(a)$ is $(-2, 0 ). $

  Let $ \vartheta_0 : \bC_+ \to \Pi_0(a)$ denote the conformal map satisfying 
  $$  \vartheta_0 (-1) = -1, \quad  \vartheta_0(+1) = 0 \quad\text{and}\quad  \vartheta_0(\infty) = \infty . $$
  We point out that the domain $\Pi_0$ is locally connected and therefore the given conformal mapping extends continuously  to the boundary \cite[p.6]{carleson-gamelin}.
  
  Taking into account the symmetry of  $ \Pi_0(a)$ we observe that
  $$  \vartheta_0(0) = -\frac{1}{2}+\frac{a}{2}i. $$
Next define  $$  -\nu(a) =   \vartheta_0^{-1}(-1+ai). $$
Reflecting $\vartheta_0$ along the half axes $ (-\infty,\nu(a)) $ respectively $ (-1 + i a,-1+i\infty) $  we obtain a conformal extension of $ \vartheta_0 $,
  $$ \vartheta : \bC \sm [ -\nu(a),\infty  ) \to 2 \Pi_0(a) $$
  satisfying
  $$ \vartheta(1) = 0 , \quad \vartheta(-\nu(a)) = -1 + i a ,  \quad  \vartheta(\infty) = \infty . $$
  
  \begin{lemma}  For any $ a > 0 ,$ the mapping
    $$ Q(z) = \vartheta^{-1}(2 \vartheta_0(z) ), $$
coincides with  the quadratic polynomial
$$ 
(\nu(a) +1) z^2 - \nu(a) .  $$
Putting  $ \xi = (\nu(a) +1)  $ and  $ \varphi(z) = z/\xi , $ 
        yields
  $$  \varphi^{-1}Q ( \varphi(z)) = z^ 2 - \nu(a)(\nu(a)+1) ,
 \quad\text{ and }\quad 
   \varphi^{-1}Q  (\varphi(\xi ))
   = \xi  .$$
   \end{lemma}
   \proof
    First observe that $ Q : \bC_+ \to  \bC \sm  [ -\nu(a),\infty  ) $ defined as
  $$ Q(z) = \vartheta^{-1}(2 \vartheta_0(z) )$$
  is a conformal mapping satisfying 
  $$ Q(1) = 1, \quad Q(0) = -\nu(a)  , \quad Q(\infty) = \infty . $$
  Indeed  we have
  $\vartheta^{-1}(2 \vartheta_0(0 ) ) 
  = 1,$
  and 
  $$ \vartheta^{-1}(2 \vartheta_0(1 ) )=
   \vartheta^{-1}\left(2 \cdot  (-\frac{1}{2}+\frac{a}{2}i )\right)  = \vartheta^{-1}(-1+a i ) =
  -\nu(a) .$$
   Hence   $ Q : \bC_+ \to  \bC \sm  [ -\nu(a),\infty  )$
%
%
   coincides with
   $$ 
   (\nu(a) +1) z^2 - \nu(a)  . $$
   Finally, we pass to a monic quadratic polynomial $P(z)$. By a direct calculation, putting  $ \xi = (\nu(a) +1)  $ and  $ \varphi(z) = z/\xi , $ 
        yields
  $$  \varphi^{-1}Q ( \varphi(z)) = z^ 2 - \nu(a)(\nu(a)+1) ,
 \quad\text{ and }\quad 
   \varphi^{-1}Q  (\varphi(\xi ))
   = \xi  .$$
 
   \endproof
 \paragraph{Remark:} Expressing  $ \lambda(a)  =   \nu(a) (\nu(a) +1) $ in terms of the fixed point $ \xi (a)   $ gives  
  $  \lambda (a)  = \xi (a) ^2 - \xi(a)  . $ In this case $\nu>1$ is equivalent to $\lambda>2$ and $\xi>2$.
 Thus  we found the dependence ($a>0$)
 $$ a \mapsto \nu(a) \mapsto \xi(a) \mapsto \lambda ( a) $$
 such that with 
  $ \varphi(z) = z/\xi(a)   $  and  $Q(z) = \vartheta^{-1}(2 \vartheta_0(z) ), $
   $$  P(z):=\varphi^{-1}Q ( \varphi(z)) = z^ 2 -  \lambda(a) 
 \quad\text{ and }\quad 
  P(\xi)= \varphi^{-1}Q  (\varphi(\xi ))
   = \xi  .$$
  Note that dynamical systems generated by $P(z)$ and $Q(z)$ are equivalent.

   Next we reverse the above chain of dependence.
   Here we use the following lemma:
    \begin{lemma}\cite{sodin-yuditski-1990} If $ a > 0 $ then  $\nu(a)  > 1 . $ 
  The mapping 
  $$a \to \nu(a) $$
  is increasing, continuous  and satisfies   $\lim_{a\to 0} \nu(a) = 1,$    $\lim_{a\to \infty} \nu(a) = \infty . $
  \end{lemma}

\begin{figure}
		\begin{center}   
\includegraphics[scale=0.6]{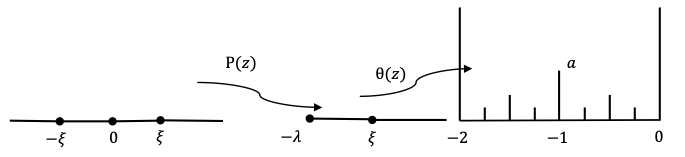}
\end{center}
\caption{Conformal mappings $\bC_+\stackrel{P}{\to}\bC\setminus[-\lambda,\infty)\stackrel{\theta}{\to} 2\cdot\Pi_0(a)$}
			\label{fig:1}
		\end{figure}

 \begin{theor} For any given  quadratic polynomial
    $$P(z) = z^ 2 - \lambda    , \quad \lambda  > 2, $$ 
    there exists  a unique $ a =  a(\lambda) >0$    and conformal mappings
$$ \theta_0( z) = \vartheta_0 \circ  \varphi(z) , \quad   \theta ( z) = \vartheta \circ  \varphi(z)  , \quad  \varphi(z) = \frac{z}{ \nu(a)  +1} ,  $$
such that
$$P(z) =  \theta^{-1}(2  ( \theta_0 (z )).$$
Moreover $\xi(a) = \nu(a)  +1   $  is the positive fixed point of $P(z) $ , i.e.,
$$P(\xi(a) ) =   \xi (a).$$
  \end{theor}
  \proof See 	Figure~\ref{fig:1}.
Given $\lambda \in ( 2, \infty) $ put
$$ \xi = \frac{ 1 + \sqrt{1 + 4 \lambda }}{2} , \quad \nu(\lambda) = \xi -1  $$
Next invoke that  there exists a unique  $ a= a(\lambda) >0 $ such that 
    $$ \vartheta^{-1}(-1+a i ) = -\nu(\lambda)  .$$  
    \paragraph{Remark:}
    Thus we  found the dependence 
$$  \lambda  \mapsto  \xi( \lambda ) \mapsto  \nu(\lambda) \mapsto a(\lambda)  $$
such that
$$ z^ 2 - \lambda =  \varphi^{-1}\vartheta^{-1}(2 \vartheta_0( \varphi(z )) $$
where  $ \varphi(z) = \frac{z}{ \xi(\lambda)} . $ Or equivalently
$$ z^ 2 - \lambda =  \theta^{-1}(2  ( \theta_0 (z )).$$
\paragraph{Remark:}
Note that the conformal map  $ \theta_0 : \bC _+ \to \Pi _0(a) $ arising in the above
Theorem satisfies
$$\theta_0 (-\xi ) = - 1 , \quad  \theta_0 (\xi ) = 0, \quad \theta_0 ( \infty  ) = \infty . $$
\paragraph{Remark:}  The identity
$$P(z) =  \theta^{-1}(2  ( \theta_0 (z )) ,  \quad z \in \bC_+$$
can be rewritten as 
\begin{equation}  \label{13mar-1}
\theta (P(z)) =  2   \theta_0 (z ) ,  \quad z \in \bC_+,\quad
\text{and}\  \theta |_{ \bC_+} = \theta_0.
\end{equation}

\paragraph{B\"ottcher's Identity.}\cite[p. 34-35]{carleson-gamelin}
 Define initially  
  $$   G (z)  =  e^{ i\pi  \theta ( z )} ,  \quad z \in  \bC_ +  . $$
By \eqref{13mar-1} we obtain  B\"ottcher's identity,
\begin{equation}  \label{13mar-2}
  G (P(z)) =   e^{ i\pi  \theta ( P(z ) ) } = e^{ i\pi 2 \theta ( z )} = G(z)^ 2 , 
  \end{equation}  
which can be extended by analyticity to the domain  $\bC  \sm [-\xi, \xi ] $ with an extension by continuity to the boundary correspondence.

\paragraph{Remark:}  A point $x\in(-\xi,\xi)$ corresponds to two different points on the boundary of the domain $\bC  \sm [-\xi, \xi ] $, which we denote by $x\pm i0$ with vicinities in the lower and upper half planes respectively.
 \begin{cor} 
Let 
 $E_0=\theta_0^{-1}([-1,0])$. Then $\bC\setminus E_0$ is the basin of attraction of infinity for the dynamical system generated by $P(z)$. Respectively $E_0$ is its Julia set. 
 \end{cor}
  \proof
  We  denote by  $P^{\circ n} $ the $n-$fold iterate of $P$.
  Note that the unit disc $\bD$ is the basin of attraction of zero for the mapping $w\mapsto w^2$.
Respectively, due to \eqref{13mar-2} we have
$$
P^{\circ n}(z)\to\infty\ \text{for all}\ z\in \bC\setminus E_0.
$$
Alternatively, since $|G(z)|=1$ if $z\in E_0$, we have $|P^{\circ n}(z)|\le \xi$.
  
     \endproof

\subsection{Green's Mapping}
Fix  $ \lambda \in \bR   $ with $ \lambda >2 . $  Let $P(z) = z^2 - \lambda. $ Let $ \xi $ denote its non-negative fixed point, $P(\xi ) = \xi . $
Let $\eta > 0 $ satisfy
$P(\eta ) = -\xi . $
\begin{enumerate}
\item
 $ \lambda >2  $  if and only if   $ \xi > 2 .$
\item 
If   $ \xi > 1 + \sqrt{2}$,  then $  \eta  > 1 .$
\item 
Let $ E_0 \sb \bR $  be  the Julia set of the quadratic polynomial $P :\bC \to \bC . $  We have,
$$  E_0 \sb [- \xi,  \xi] \setminus (-\eta, \eta) .$$
\end{enumerate}

    \paragraph{Green's Fundamental domain.}  The range of $   e^{ i\pi  \theta ( z )} ,   z \in  \bC   \sm [-\xi, \xi ]  $ forms Green's Fundamental domain.
    We put 
   $$     H = \{  e^{ i\pi  \theta ( z )} :  z \in  \bC   \sm [-\xi, \xi ] \} $$

   \paragraph{Green's mapping.}
   We define Green's mapping   $ T : H \to \bC  \sm [-\xi, \xi ] $  to be the inverse of $G : \bC   \sm [-\xi, \xi ] \to  H . $  Thus
   $T$ is defined by $ T( G (z)) = z $, explicitly in terms of $ \theta $  
    $$  T( e^{ i\pi  \theta ( z )}) = z  ,  \quad z \in  \bC   \sm [-\xi, \xi ]  . $$

    \begin{prop}  Given   $ \lambda \in \bR   $ with $ \lambda >2  .$ Let $E_0$ denote the Julia set of the quadratic polynomial  $P(z) = z^2 - \lambda. $
  There exists a unique $ a > 0 $ such that
  \begin{enumerate}
    \item 
  $\theta  (\lambda ) =   -1  + i a   $
  \item 
  $ \theta^{-1}\{ (-1, 0) \}  = E_0 $, or, equivalently, $T(\bT ) = E_0 $.  
   \end{enumerate}  
  \end{prop}

  \begin{prop} Let $P(z) = z ^ 2 -\lambda . $ Let $G$ denote the  B\"ottcher function of P, and let $E_0 $ denote the Julia set of $P$.
    Then $|G(z)| = 1 $ for $z\in E_0$ and   $|G(z)| <  1 $ for $z\notin E_0$.
    \end{prop}
    \paragraph{Remark:} In view of the above proposition,
    $$- \log |G(z) | =\pi\Im\theta(z)$$
    is just the Green's function of $\bC \sm E_0 . $ 
   
\subsection{The Poincare Function }
Fix  $ \lambda \in \bR   $ with $ \lambda >2 . $  Let $P(z) = z^2 - \lambda. $ Let $ \xi$ denote its non-negative fixed point. Then
$$ P(\xi ) = \xi , \quad  P^{\prime} ( \xi ) = 2 \xi . $$
Next put $\varphi(z)=z+\xi$ and
$$\hat P (z) = \varphi^{-1}(P(\varphi(z) ) ). $$
Then,  by direct calculation,
$$ \hat P (z) =  z ( z+\rho ), \quad \hat P(0) = 0 , \quad \hat P^{\prime} (0) = 2 \xi . $$
Below we employ the abbreviation  $ \rho = 2 \xi .$ 
Note that since  $ \xi >  2 $  we have $ \rho = 2 \xi  > 4 .$  

Thus the dynamics of the polynomial $P(z) $ is equivalent to the one of  $ \hat P (z) $ whose fixed point  equals $0$ and whose multiplier  $\rho$ coincides with the multiplier of $P$ at its fixed point $\xi$.    

\begin{theor} [see e.g. \cite{derfl, grabner}] There exists an  entire  function $ F$ which obeys the Poincare equation 
  $$ \hat P(F(z) ) = F(\rho z ),
\quad F(0) = 0 , \quad F^\prime (0 ) =1.
  $$
\end{theor}
Thus  Poincare's function $F$ conjugates $ \hat P $ to multiplication by $\rho ,  $  where   $ \rho = 2 \xi $   and  $ \xi$ is the non-negative fixed point of $P . $  We point out that
Poincare's function $F$ is explicitely given by
$$ F(z) = \lim_{n \to \infty}  \hat P^{\circ n } ( z/\rho^n) . $$
In our present context the above formula for Poincare's function follow at once from its comb-domain representation, presented below. 

\subsubsection{The Comb Domain Representation of $F$ }
\begin{figure}
\begin{center}
\includegraphics[scale=0.6]{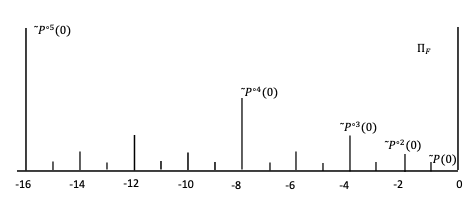}
\end{center}
\caption{Comb domain $\Pi_F$}
  \label{fig:4}
		\end{figure}

Define  the sequence $ h(-l) , l \in \bN $ by decreeing  that
 for any odd(!)  integer $ k $, and $ m \in \bN ,$ 
\begin{equation}\label{eq19-1}
- \xi \cosh (\pi h(-k) ) = P(0) \quad\text{and}\quad \xi \cosh (\pi h(-k 2^ m ) ) = P^{\circ (m+1)}   (0) .
\end{equation}
See Figure \ref{fig:4}.
We call $ \{ \psi +i \phi : \psi \in (-\infty , 0 ) , \phi > 0 \} , $ the   
 left upper quarter of the complex plane.
Let $\Pi_F$ be the  comb domain over the base interval $ (-\infty , 0 )  $ obtained by removing
from the left upper quarter of the complex plane
the following sequence of  intervals
$$ [ -l , -l +i h(-l)  ] , \quad l\in \bN .  $$
Define $ \Psi_F : \bC_+ \to \Pi_F $ to be the conformal (= analytic and bijective) map satisfying 
$ \Psi_F (0) = 0 , \quad \Psi_F(\infty) = \infty , $ and
$$  \frac{d}{dw} ( \xi \cos ( \pi \Psi (w) ))_{|_{w = 0 }} = 1 . $$
\paragraph{Remark: }
We note in passing that for $ w \in \bC_+$ with $|w|$ small enough,
$$  \Psi_F (w)  = \tau  w^{1/2} + \text{  higher order powers of $ w $ , } \quad\text{where}\ \tau =  \frac{2i}{\sqrt{{\pi \xi} }}.  $$

\begin{figure}
\begin{center}
\includegraphics[scale=0.6]{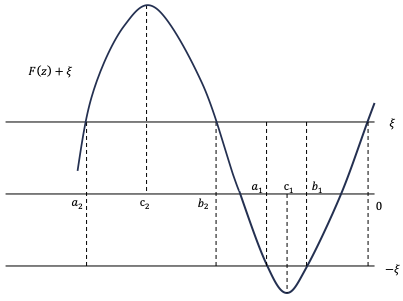}
\end{center}
\caption{Graph of  the Poincare function $F(z)$ shifted by $\xi$}
  \label{fig:3}
		\end{figure}

\begin{theor} The Poincare function  $ F: \bC_ + \to \bC $
  satisfying  $ F(0) = 0 , $ $F(\infty ) = \infty $ and   $ F^\prime (0 ) =1 $ admits the
  following representation in terms of the conformal map   $ \Psi_F : \bC_+ \to \Pi_F $
  onto the comb domain $\Pi_F$,
  $$ F(w ) +\xi =  \xi \cos ( \pi \Psi_F (w) )), $$
 (see Figure \ref{fig:3}) where $ \Psi_F (0) = 0 , \quad \Psi_F(\infty) = \infty , $ and
$$  \frac{d}{dw} ( \xi \cos ( \pi \Psi_F (w) ))_{|_{w = 0 }} = 1 . $$
\end{theor}
\begin{proof}
Let  $\cC_{k}=\{c^{k}_{j}: (P^{\circ k})'(c^{k}_j)=0\}$ be the collection of the critical points of $P^{\circ k}$.
Since $P^{\circ (k+1)}(z)=P^{\circ k}(P(z))$, we have
$$
(P^{\circ (k+1)})'(z)=(P^{\circ k})'(P(z))P'(z)
$$
That is, $\cC_{k+1}=\{0\}\cup\{  c^{k+1}_{j}: P(c^{k+1}_j)=c^k_j\}$. Respective critical values are
$$
P^{\circ (k+1)}(0)\quad\text{and}\quad P^{\circ k}(c^k_j).
$$
Since $P(z)$ increases on the positive half axis and decreases on the negative half axis, in passing from $k$ to $k+1$, the critical values of the previous generation  are listed in the same order
on the positive half axis and in the opposite order on  the negative  one. In summary $\hat P^{\circ n}(z/\rho^n)$ can be restored in terms of its critical values listed in a proper oder
 as
$$
\hat P^{\circ n}(z/\rho^n)+\xi=\xi\cos\pi \Psi_{\hat P_n}(z),
$$
where $\Psi_{\hat P_n}$ is the conformal mapping from $\bC_+$ to the domain $\Pi_{\hat P_n}$,  see e.g. \cite{sodin-yuditski-92, eremenko-yuditski-2012}, with the following normalization
$$
\Psi_{\hat P_n}(0)=0,\quad \Psi_{\hat P_n}(\infty)=\infty,\quad (\xi\cos\pi \Psi_{\hat P_n})'(0)=1.
$$
Here $\Pi_{\hat P_n}$ is the  comb domain over the base interval $ (-2^n , 0 )  $ obtained by removing
from the half strip, i.e., 
$ \{ \psi +i \phi : \psi \in (-2^n , 0 ) , \phi > 0 \} , $
the following sequence of  intervals
$$ [ -l , -l +i h(-l)  ] , \quad 1\le  l\le 2^n-1  $$
and $h(-l)$ are given by \eqref{eq19-1}.

Evidently combs  $\Pi_{\hat P_n}$ converges to  $\Pi_{F}$ in the Carath\' eodory sense.
\end{proof}

\paragraph{Defining $ f_n : \bC _+\to \cD_n$.}
Define  $c_{n   }$ by the conditions that
$$
\Psi_F (c_{n   } ) =\text { tip of the cut of $\Pi_F$ at $-n$}. $$
Note that $c_n$-s are critical points of $F(z)$, moreover  $c_{2n} $ is a   critical points of $F_{|\bR}$
corresponding to a  local maximum of $F_{|\bR}$
and $c_{2n +1} $ is a critical point of $F_{|\bR}$ corresponding to a local minimum  of $F_{|\bR}$. 
Let
$$\gamma_n^F = \Psi^{-1} ([ ih(-2n) -2n ,  i\infty -2n ]) $$
and let $\cD_n \sb \bC_+ $ be the domain defined by its boundary
$$ \pa \cD_n = \gamma_{n+1}^F \cup [c_{2n +2} , c_{2n} ] \cup \gamma_n^F  . $$
Define  $a_{n   }, b_{n   } \in \bR $ (in terms of prime ends) by the conditions that
\begin{equation*}
\Psi_F( a_{n }) = -n -0     ,\quad  
\Psi_F( b_{n }) = -n +0 .
\end{equation*}
We remark  $a_{n   }, b_{n   }, c_n  \in \bR $ are ordered  as follows
 $$ a_{2n +2}< c_{2n +2} < b_{2n +2}< a_{2n+1}< c_{2n +1} <  b_{2n +1} <a_{2n}<c_{2n   }<  b_{2n} $$
and that 
 $$F( a_{2n }) = 0   ,\quad  F( b_{2n }) = 0 , \quad  F(c_{2n+1}) + \xi = -\lambda  . $$ 
We define   $ f_n : \bC_+ \to \cD_n$  to be the conformal map  satisfying
\begin{equation}\label{eq17mar-5}
 f_n(-\xi)  = b_{2n+2} \quad f_n (\xi ) = a_{2n} , \quad\text{and} \quad f_n(\infty) = \infty . 
 \end{equation}
By symmetry we have $  f_n(0)  = c_{2n+1} .$ 
Summing up, 
$$F(\cdot) + \xi : \cD_n \to \bC \setminus [-\lambda , \infty) $$
and
$$F \circ f_n  + \xi : \bC_+ \to \bC \setminus [-\lambda , \infty) $$
are conformal
maps, and $F \circ f_n (0) + \xi = - \lambda $  , $ F \circ f_n (\xi )   + \xi  = \xi . $
Separately  $ P : \bC_+ \to \bC \setminus [-\lambda , \infty)$ is conformal
satisfying $P(0) = - \lambda, P(\xi) = \xi . $
Hence
$$ F(f_n(z) ) + \xi = P(z) , \quad z \in \bC_+.$$
\begin{figure}
\begin{center}
\includegraphics[scale=0.26]{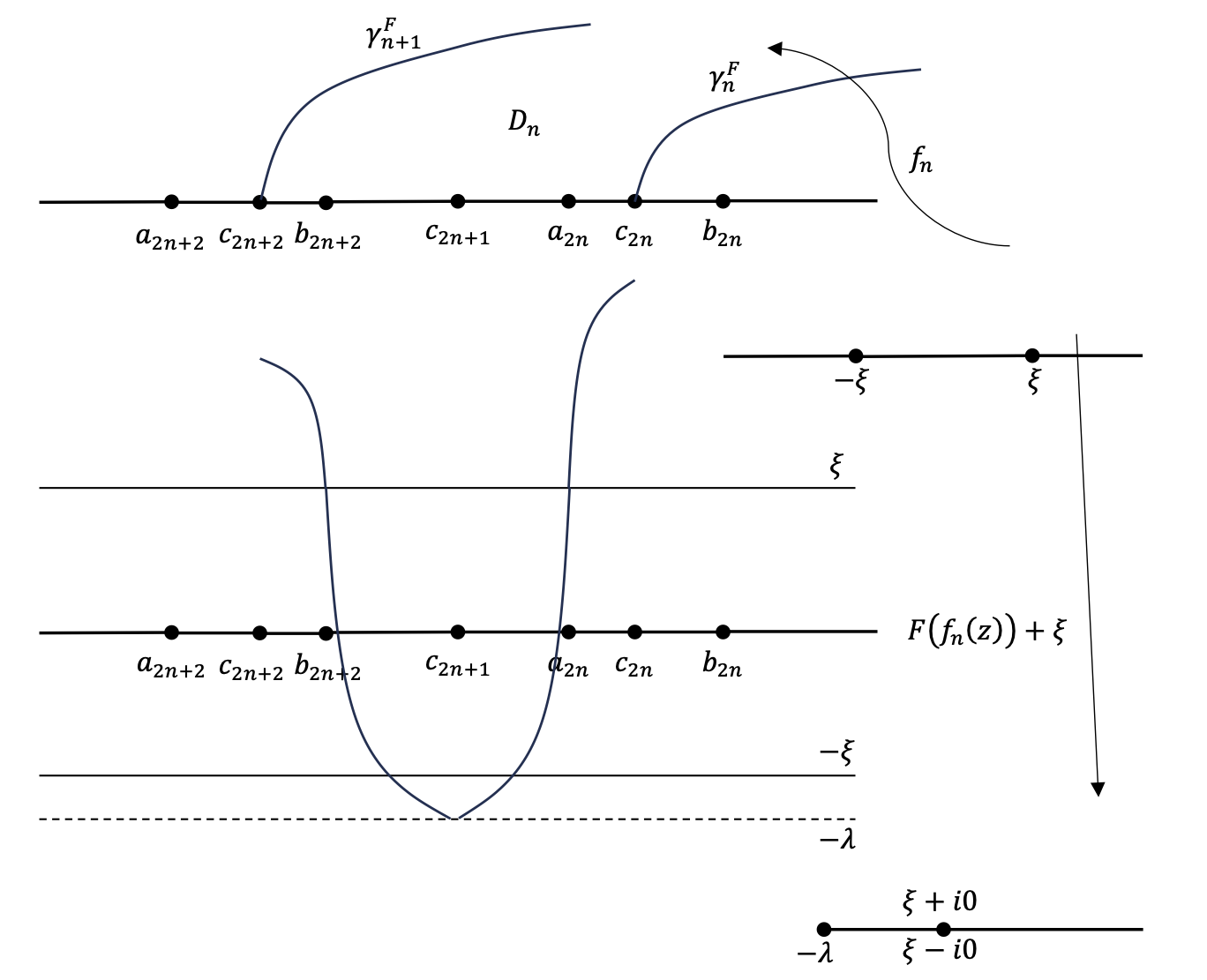}
\end{center}
\caption{$f_n:\bC_+\to \cD_n$, graph of $F(z)+\xi$ over  $[c_{2n+2},c_{2n}]$ and their superposition}
  \label{fig:5}
		\end{figure}  

\subsection{Defining $\Delta  : \bC_+ \to  \Pi $ in terms of $F$ and $\theta . $}\label{2.5}
\begin{lemma} The composed mapping $ \theta( F(\cdot) + \xi): \cD_n \to  2\cdot \Pi_0(a) $ is
  conformal     
  and
  $$ 2\theta( F(z) + \xi)  =  \theta( F(\rho z) + \xi) $$
  \end{lemma}
  \proof  Combine the fact that
  $ \theta : \bC \sm ( -\lambda , \infty) \to 2\cdot\Pi_0(a) $ and
  $F(\cdot) + \xi : \cD_n \to \bC \setminus (-\lambda , \infty) $ conformal
    maps. Taking into account the defining properties  of the B\"ottcher and Poincare maps gives 
  $$ 2\theta( F(z) + \xi)  =    \theta( P(F( z) + \xi)) =  \theta( \hat P(F( z) ) + \xi)) . $$
Invoking the  Poincare equation  we find 
  $$  \theta( \hat P(F( z) ) + \xi)) =  \theta (F( \rho z) ) + \xi)  . $$ 
  \endproof

  Next  the following sequence of conformal maps  $\Delta_{|\cD_n}  : \cD_n \to 2\cdot \Pi_0(a) -2n $ defined  by
  $$ \Delta_{|\cD_n} (z) = \theta( F(z) + \xi)-2n \qquad z \in  \cD_n ,  $$
  gives rise to
  $$ \Delta :  \bigcup_{n = 0 } ^\infty  \cD_n \to \bigcup_{n = 0 }^\infty  \{2\cdot \Pi_0(a) -2n\}   $$
  which by continuity extends to a conformal 
  map
  $ \Delta :  \bC_+  \to  \Pi    .$
  \endproof
\begin{lemma} The conformal map  $ \Delta  : \bC_+ \to  \Pi $  satisfies 
  $$ 2 \Delta (z)  =  \Delta (\rho z) , \qquad z \in \bC_+ . $$
    \end{lemma}
    \proof
    Since $ 2  \Pi = \Pi , $ we have  $2 \Delta (z)  \in \Pi . $  Hence forming
  $ \Delta^{-1} (  2 \Delta (\cdot) ) :  \bC_+ \to  \bC_+ , $
  yields a well defined conformal
    map satisfying
$$  \Delta^{-1} (  2 \Delta (\bR) )  = \bR , \quad 
\Delta^{-1} (  2 \Delta (0) ) = 0 , \quad  \Delta^{-1} (  2 \Delta (\infty) ) = \infty  . $$
Consequently there exists $ \rho _0  \in \bR $ such that
$$ \Delta^{-1} (  2 \Delta (z) ) = \rho_0 z . $$
In view of the previous Lemma we find  $ \rho _0 =  \rho  . $ Indeed for $ w \in \cD_0$
\begin{equation}
  2\Delta(w) =  2 \theta (F(w) + \xi ) = \theta ( P(F(w) + \xi) ) =  \theta ( F(\rho w) + \xi ) =  \Delta( \rho w),
\end{equation}
implying  $ \rho _0 =  \rho  . $
\endproof

\begin{lemma}
The conformal mapping 
$$ f_n =  (F(\cdot) + \xi) ^{-1} \circ P :  \bC_+  \to \cD_n , $$
satisfies,
$$
f_n( z) = \Delta^{-1} (\theta (P(z)) -2n)   , \qquad z \in \bC_+ .  $$

  \end{lemma}
\section{Preparation I I}
  \subsection{The Good Set $F_N . $}
%
%

  Let $ N \in \bN . $ Let $ \psi \in [0,1]$ with dyadic expansion   $\psi = \sum_{k = 1}^\infty \varepsilon_k 2^{-k},$
  where $ \varepsilon_k \in \{0, +1\}. $  We say that $ x \in F_N $ if $\varepsilon = ( \varepsilon_k) $ satisfies the following conditions:
  For every $ k \in \bN $ the hypothesis
  $\varepsilon_k =1 ,\varepsilon_{k+1} = 1 , \dots , \varepsilon_{k+N} = 1 $ implies that
  $ \varepsilon_{k+N+1} = 0 ;  $
  and    $\varepsilon_k =0 ,\varepsilon_{k+1} = 0 , \dots , \varepsilon_{k+N} = 0 $ implies that $ \varepsilon_{k+N+1} = 1 . $

  Note  that  $F_N  $ is a closed subset of $[0,1]$ without isolated points. In this section we prove that
  \begin{prop} \label{p19-6}
    There exists  $ c> 0 $  such that for any $N\in \bN $ 
    $${\rm Hdim}(F_N) > 1 - c/N . $$ 
\end{prop} 
\subsubsection{Representing $F_N$}  

Let $ m \in \bN $ and $\alpha \in \{0,1\}^m . $ We let $I$ denote the unique dyadic interval of length $ 2^{-m} $ whose left endpoint equals
  $\alpha_1 2^{-1} + \cdots + \alpha_m 2^{-m} .$
  Thus 
  $$ I = \alpha_1 2^{-1} + \cdots + \alpha_m 2^{-m} + \sum_{k = 1}^\infty \varepsilon_{m+k} 2^{-m-k} :  \varepsilon_{m+k} \in  \{0, +1\} .$$
  We say $\alpha$ is the index of the given interval $I$, in short $ind(I)=\alpha$.
  In what follows $I_0$ and $I_1$ are the intervals with indices $(\alpha,0)$ and $(\alpha,1)$ respectively.
  
  Given any collection of dyadic intervals $\mathcal  B$  set
  $$
  \mathcal  B^*=\cup_{J\in \mathcal B}J.
  $$
  
  \paragraph{Defining $\mathcal E^2 (I),\mathcal K^2(I)$. }
  We say that  $\alpha \in \{0,1\}^m  $ is of Type $0$ if $\alpha_m  = 0. $
 For a  Type $0$ interval  $I $, $ind(I)=\alpha$, define
 \begin{equation}
 \mathcal E^2 (I) = \{ I_{( 0,0)} ,   I_{( 1,1,1)} \} \qquad \mathcal K^2 (I) = \{ I_{(0,1)} ,   I_{(1,0)} , I_{(1,1,0)} \} 
\end{equation}
  We say that  $\alpha \in \{0,1\}^m  $ is of Type $1$ if $\alpha_m  = 1. $
 For a  Type $1$ interval  $I $ define
  \begin{equation}
 \mathcal E^2 (I) = \{ I_{( 1,1)} ,   I_{( 0,0,0)} \} \qquad \mathcal K^2 (I) = \{ I_{( 0,1)} ,   I_{(1,0)} , I_{(0,0,1)} \} 
\end{equation}
Taking unions we put
$$  K^2 (I) = \bigcup_{J \in \mathcal  K^2 (I)} J ,\qquad  E^2 (I) =\bigcup_{J \in \mathcal  E^2 (I)} J , $$
satisfying 
$$ K^2 (I) \cap   E^2 (I)  = \es , \qquad  K^2 (I) \cup   E^2 (I)  =I , $$

$$|  K^2 (I) | \ge \frac58 |I| , $$
and 
$$ |J| \le \frac 14 | I| , \quad \text{for} \quad J \in \mathcal  K^2 (I) . $$

\paragraph{Defining  $\mathcal E^{N+1} (J_*),\mathcal K^{N+1}(J_*)$. } 
Suppose the collections of intervals $\mathcal E^N(J_*)$ and   $\mathcal K^N(J_*)$  are defined for any dyadic interval $J_*$.

\paragraph{Part 1.}
Fix $J \in  \mathcal E^N(J_*).$ If $J$ is of Type 1 then we decree
$$ \cF (J) =\{J_1\} , \qquad \cG (J)=\{J_0\} . $$
If  $J$ is of Type 0 then 
$$ \cF (J) =\{J_0\} , \qquad \cG (J)=\{J_1\} . $$
Next form
$$\mathcal E^{N+1} _{ \mathcal E^{N}}(J_*) = \{ \cF (J) : J \in  \mathcal E^N(J_*)\}, \qquad \mathcal K^{N+1} _{ \mathcal E^{N}}(J_*) = \{ \cG (J) : J \in  \mathcal E^N(J_*)\} . $$

\paragraph{Part 2.} Fix  $J \in  \mathcal K^N(J_*)$, such that the index of $J$ ends with $(1,0)$. Then form the indices
$$\alpha_k = (0, \dots , 0)  = \{0\}^k  $$
where $k \le N +1. $ We put
$$\cG(J) = \{ J_1, J_{(\alpha_1,1)},  J_{(\alpha_2,1)}, \dots ,  J_{(\alpha_N,1)} \} , \qquad  \cF(J) = \{ J_{\alpha_{N+1}}\}. $$
Next consider  $J \in  \mathcal K^N(J_*)$, such that the index of $J$ ends with $(0,1)$.
 Then let  $\beta_k = (1, \dots , 1)  = \{1\}^k  $ where $k \le N +1. $
We put
$$\cG(J) = \{ J_0, J_{(\beta_1,0)},  J_{(\beta_2,0)}, \dots ,  J_{(\beta_N,0)} \} , \qquad
\cF(J) = \{ J_{\beta_{N+1}}\}. $$
Clearly we have $|I| \le \frac12 |J| $ for $I \in \cG(J)$ and  $ |\cF (J)| = 2^{-n-1} |J| ,$
Next form
$$\mathcal E^{N+1} _{ \mathcal K^{N}}(J_*) = \{ \cF (J) : J \in  \mathcal K^N(J_*)\}, \qquad \mathcal K^{N+1} _{ \mathcal K^{N}}(J_*) = \{ \cG (J) : J \in  \mathcal K^N(J_*)\} . $$
\paragraph{Part 3.} Take unions 
$$\mathcal E^{N+1} (J_*) = \mathcal E^{N+1} _{ \mathcal K^{N}}(J_*)  \cup \mathcal E^{N+1} _{ \mathcal E^{N}}(J_*)$$
and 
$$\mathcal K^{N+1} (J_*) = \mathcal K^{N+1} _{ \mathcal K^{N}}(J_*)  \cup \mathcal K^{N+1} _{ \mathcal E^{N}}(J_*)$$

\paragraph{ Part 4. Backtracing:}
The above definition yields that 
$$|\mathcal E^{N+1} (J_*)^* |= \frac12 |\mathcal E^{N} (J_*)^*| +\frac{1}{2^{n+1}} | \mathcal K^{N} (J_*)^*| \le \frac12 |\mathcal E^{N} (J_*)^*| +\frac{1}{2^{N+1}} |J_*| $$
and 
$$ |J | \le \frac{1}{2^{N+1}} |J_*|, \quad 
\text{for} \quad J\in \mathcal E^{N+1} (J_*) \cup \mathcal K^{N+1} (J_*) . $$
Since $|\mathcal E^{2} (J)^* | \le \frac38 |J|$ for any dyadic interval $J$ it follows from the above recursion that
$|\mathcal E^{N+1} (J_*)^* | \le \frac12  |J_*| $. 
\endproof

Recall $J$ is a half open interval,  by $\bar J$ we mean its closure.

  \begin{lemma}
    There exist finite collections of pairwise dyadic intervals
    $\{ \mathcal K _ k^N , k \in \bN \} $ such that with $ K_ k^N = \bigcup_{J \in    \mathcal K _k ^N } \overline{J} $
    the following conditions hold true
\begin{enumerate}
\item  $ F_N \sbe K_k^N  $ for $ k \in \bN $.
\item  $  F_N = \bigcap_{ k \in \bN } K_k^N . $ 
\item For $J \in    \mathcal K_k^N$,  $| J \cap   K_{k+1}^N| \ge \frac12 |J| . $
\item If  $J \in    \mathcal K_k^N$, $I \in    \mathcal K_{k+1}^N$ and  $ I \cap J \neq \es $ then  $I \sb J $ and $|I| < 2^ {-N } |J| .$
\end{enumerate}
\end{lemma}  
\proof
Start by putting  $ \mathcal K_0^N = \{ [0, 1/2 ], [1/2, 1] \} .$ In the first step we define
$$
 \mathcal K_1^N = \{  \mathcal K^N ( I ) :  I \in  \mathcal K_0^N \}  \qquad  \mathcal E_1^N = \{  \mathcal E^N ( I ) :  I \in  \mathcal K_0^N \} 
 $$
 Having defined $\mathcal K_m^N  ,    1 \le m \le k $  we put
 $$
 \mathcal K_{k+1}^N = \{  \mathcal K^N ( I ) :  I \in  \mathcal K_k^N \}  \qquad  \mathcal E_{k+1}^N = \{  \mathcal E^N ( I ) :  I \in  \mathcal K_k^N \} 
 $$
 In view of the above  we obtain  $| J \cap   K_{k+1}^N| \ge \frac12 |J| , $
   $J \in    \mathcal K_k^N$, 
   and if   $I \in    \mathcal K_{k+1}^N$  with   $ I \cap J \neq \es $  we obtain   $I \sb J $ and $|I| < 2^ {-N } |J| .$
   
Next  we have   $  K_{k+1}^N  \sb  K_{k}^N  ,$ and $ E_{k+1}^N  \sbe  K_{k}^N . $ Moreover the sets  $ E_{k}^N , k \in \bN , $  are pairwise disjoint and
 $$ \bigcup_{k \in \bN }  E_{k}^N \cap F_N = \es . $$
 Consequently $ K _{k}^N  \sp  F^N $  for  $k \in \bN , $  and
 $$ \bigcap_{k \in \bN }  K_{k}^N  = F_N  . $$

 \endproof

  \subsubsection{Hausdorff measure estimates}
  \begin{prop}
    There exists  $ c> 0 $  such that for any $N\in \bN $ 
    $$ \Lambda_{1- \frac{c}{N}} (F_N) >  0  $$ 
  \end{prop} 
  
   \proof
  Use the criterion in Pommerenke's book \cite[p.226] {pommerenke} in combination with the representation of $F_N$ given above.
  Define 
  $$ \mathcal A_k = \{ \overline{J \sm {\mathcal E ^ N (J)}^* } : J \in \mathcal K_k^N \} ,  $$
  and set $A_k = \bigcup_{ L \in  \mathcal B_k } L . $ Note that $ \mathcal A_k $ is a finite collection of pairwise disjoint closed sets.  Moreover if
  $ L_k \in  \mathcal A_k$  and $ L_{k+1} \in   \mathcal A_{k+1}$   satisfy
  $ L_k \cap L_{k+1} \ne \es $ then  $L_{k+1} \sb L_k$ and
  $${\rm diam }(L_k)
  \le 2^{-N}{\rm diam }(L_{k+1}) .$$ 
  If $ L \in  \mathcal A_k $ then  $$| A_k \cap L | \ge \frac12 |L| $$
  and
  $$ \bigcap _k A_k = F_N . $$
  Now Theorem 10.5 in  Pommerenke's book \cite[p.226] {pommerenke} implies 
  $$\Lambda_\alpha (F^N ) > 0, \quad \text{for} \quad  \alpha = \frac{\log( 2^{N-1})}{\log( 2^{N})} = 1-\frac{1}{N}. $$
  
  \endproof
  
 \subsection{A  compact $K \sb \bC_+$ defined by $F_N$.}\label{sub18mar-1}
\begin{figure}
 \begin{center}
\includegraphics[scale=0.4]{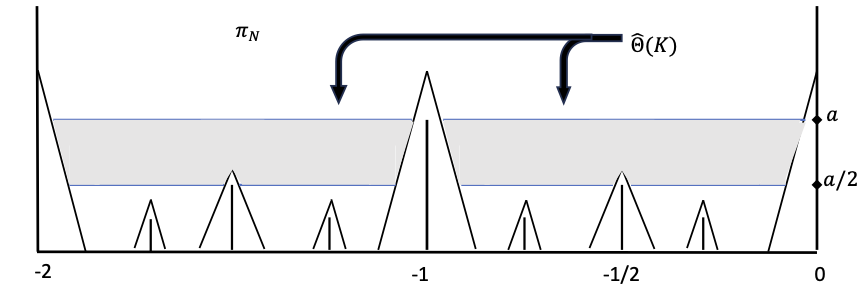}
\end{center}
 \caption{Image of the compact $K$ in the domain $2\Pi_0(a)$}
			\label{fig:7}
		\end{figure}

Given  $\psi \in F_N $  and $ 0 < \gamma  < 1/2 $ set
  $$
  V_{\psi , \gamma} = \{ x + iy \in (-1,0)\times \bR_+: | x + \psi | < \gamma y \} .$$
  Note that  $V_{\psi , \gamma}$ is a truncated cone, with vertex  at $-\psi$, contained in  the half strip $\{(-1,0)+i\bR_+\}$. Its opening angle equals
  $2\arctan (\gamma)$ and its axis of symmetry parallel to the $y-$axis.
  \begin{lemma} If $ \gamma < ( a 2^{N+1})^{-1} $ and  $\psi \in F_N $  then   $ V_{\psi , \gamma} \sbe \Pi_0(a) . $
  \end{lemma}
  Let $ \gamma < ( a 2^{N+1})^{-1} $. We form
  $$\pi_{0,N} = \cup_{\psi \in F_N} V_{\psi , \gamma} $$ 
  and let $\pi_{N}$ be the union of  $\pi_{0,N} $ and its reflection at the line $\{-1 +i y :  y \in \bR_+\}  $. Next  put
  $$ B = B(N)= \pi_{N} \cap \left\{\bR + i[{a}/{2} , a ]\right\} \subset 2\cdot \Pi_0(a)
  $$
  See Figure~\ref{fig:7}. Finally we put
  $$K_0=K_0(N) =  \hat \theta ^{-1}(B(N)),\quad K(N)=\overline{K_0(N)}, $$
  where  $\hat \theta (z) = \theta (P(z)) , z \in \bC_+ .$
  
  \begin{lemma}\label{l21-1}
    \begin{enumerate}
    \item 
    $K  =  \overline{\hat \theta ^{-1}(B)}  $ is a compact subset of $ \bC_+$.  
    \item Let $L_{\psi}=\hat\theta^{-1}([-\psi+ia/2, -\psi+ia])$. Then $L_\psi \subset K(N)$ for all $\psi\in F_N$.
  \end{enumerate}
  \end{lemma}
\proof
$F_N$ is a Cantor type set (closed without isolated points). We have 
$$
(-1,0)\setminus \{-\psi: \psi\in F_N\}=\cup_{j} I_j
$$ 
with a system of complementary open intervals $I_j$. Therefore $\pi_{0,N}$ represents a complement to the union of symmetric triangles $\delta_j$  based on $I_j$ and having the 
same angle at their bases (depending on $\gamma_N$), see Figure \ref{fig:7}. Since $|I_j|\to 0$, as $j\to\infty$, only a finite number of $\delta_j's$ intersect the strip 
$ \left\{\bR + i[{a}/{2} , a ]\right\}$. Thus $B=B(N)$ is a pre-compact subset of  $2\cdot\Pi_0(a)$. Respectively its preimage $\hat \theta ^{-1}(B(N))$ is pre-compact in $\bC_+$. Note that the number of connected components of $K(N)$ is growing with $N$, but it is always finite.

Since each interval $[-\psi+ia/2, -\psi+ia]$, $\psi\in F_N$, is contained in $B(N)$, we get $L_\psi\in K(N)$ for all $\psi\in F_N$.
\endproof

\section{Integral Estimates}
  In this section by several steps of reduction  we prove our main Theorem \ref{main}.
    \subsection{First Reduction}
    Fix $\psi \in (0, 1 )  $.  We  consider the radial path $ R_{\psi}$ defined by 
    $$ R_{\psi} (h) =  e^{ i\pi ( -\psi +ih) }, \quad h < a, $$ in $H$, the vertical path  $$ s_{\psi} (h) =   ( -\psi +ih), \quad h < a,   $$
    in  $ \Pi_0(a)$
    and we let $$ \gamma_\psi (h) =  \theta^ {-1} (s_{\psi} (h)), \quad h < a  .$$  
    Thus $ \gamma_\psi $  is a path in $ \bC_ + , $   defined as the image of  $s_{\psi} $ under the inverse map $ \theta^ {-1} $.  We have
    $$  R_{\psi} = 
    e^{ i\pi  s_\psi }  =  e^{ i\pi  \theta ( \gamma_\psi )} .$$

 \paragraph{Remark:} The condition
 $$
 I_T:=\int _{ R_{\psi}} |T^{''} (w) | \cdot |d w |<\infty
 $$
 implies evidently that $T'(w)$ is bounded on $R_\psi$. Indeed, for $r\in (e^{-a},1)$
 $$
 |T'(r e^{-\pi i\psi})|\le \left|\int_{e^{-a\pi}}^r T''(te^{-i\pi\psi}) dt
 \right|+ |T'(e^{-\pi (a+i\psi)})|\le I_T+ |T'(e^{-\pi(a+ i\psi)})|.
 $$
 Similar estimations will be used in what follows.
    
    \begin{theor}\label{t19-5}  For any  $\psi \in (0, 1]  $, the condition 

    $$\int_ {\gamma_\psi} \frac { | \theta ^{''} (z) | }{|\theta^{ '} (z)^ 2|} \cdot  |dz|  < \infty,  $$ 
 implies that
      $$ \int _{ R_{\psi}} |T^{''} (w) | \cdot |d w |  < \infty .$$
   \end{theor}
   \proof
As before  $G(z) = e ^{\pi i\theta (z)} $. Since $T(G(z)) = z $  we have
   $$ T^\prime(G(z)) G^\prime (z) =1  $$
   and
   $$  T^{\prime \prime} (G(z)) G^{\prime 2} (z) + T^\prime(G(z)) G^{\prime\prime} (z) = 0 $$
   Note that $ T^\prime(G(z))  = 1/ G^\prime (z)  $ implying
   $$ T^{\prime \prime} (G(z)) G^{\prime 2} (z) + \frac{ G^{\prime\prime} (z)}{ G^\prime (z)} = 0 $$
   or 
   $$ T^{\prime \prime} (G(z)) G^{\prime 3} (z) + G^{\prime\prime} (z) = 0 $$
   and
   $$ T^{\prime \prime} (G(z)) =- \frac {G^{\prime\prime} (z)}{ G^{\prime 3} (z)} .$$
In its turn, we have
$G^{\prime } (z) =  i \pi e ^{i\pi \theta (z)} \theta^\prime (z) $
and
$$G^{\prime\prime } (z) =  - \pi^2 e ^{i\pi \theta (z)} \theta^{\prime 2} (z) +  i \pi e ^{i\pi \theta (z)} \theta^{\prime \prime} (z) .$$
Since $g( {\gamma_\psi} ) =   
{ R_{\psi}}$, the transformation  $w = G(z) $, $|dw| = |G^{\prime}(z)| |dz| $ gives
$$  \int _{ R_{\psi}} |T^{''} (w) | \cdot |d w |  =    \int_ {\gamma_\psi}  \frac {|G^{\prime\prime} (z)|}{ |G^{\prime 3} (z)|} |G^{\prime}(z)| |dz| = \int_ {\gamma_\psi}
\frac {|G^{\prime\prime} (z)|}{ |G^{\prime 2} (z)|}  |dz|  . $$
By the calculations at the beginning of the proof  we have
$$
\frac{G^{\prime\prime} (z)}{ G^{\prime 2} (z)}  = \left[ 1 +\frac{1}{\pi i}\frac {\theta^{\prime\prime} (z)}{ \theta^{\prime 2} (z)} \right]\frac{1}{e ^{i\pi \theta (z)}}, 
  $$
  which concludes the proof (recall $\theta(z)$ is continuous up to the boundary).
   \endproof
%
   \subsection{Second Reduction}
   Let $\Pi \sbe \bC _+$ be a comb domain over the base interval $(-\infty, 0 )$ uniquely determined by the following two conditions  
   \begin{enumerate}
   \item  $2 \Pi = \Pi$ 
   \item $\Pi$ restricted to the base interval $ (-1,0) $ coincides with $\Pi_0(a)$
   \end{enumerate}
 Recall $ \Delta :  \bC _+ \to \Pi$  denote the conformal map satisfying  $ \Delta (0) = 0 $, $ \Delta (\infty) = \infty , $  and
   $$ \Delta (  \rho w ) = 2  \Delta ( w ) , \quad w \in  \bC ^+ ,\quad \rho=2\xi. $$
   Define
   $$ f(z) = \Delta^{-1} (\theta (z) ) , \quad z\in  \bC ^+ . $$
   Put  $ \Gamma_\psi = f (\gamma_\psi ) . $ Note that then 
   \begin{equation}\label{eq17mar-3}  \Gamma_\psi = \{w \in \bC^+ :  \Delta(w) = - \psi +ih : h \in (0, a) \} .
   \end{equation} 
   
 \begin{lemma}\label{l18mar-1}  Let $\Omega_2 $ denote the domain of holomorphy of  $f . $ 
      The path  $ \gamma_\psi $ forms a compact subset in  $\Omega_2 $. 
    \end{lemma}
    \proof
    Since $f:\bC_+\to\bC_+$ and it is real on $(-\lambda,\infty)$, $f$ has an
  analytic extension in  $\Omega_2= \bC \sm (-\infty , -\lambda ]$. The path $ \gamma_\psi $ is compact in this domain.
    \endproof

            \begin{theor}\label{t19-4}  Let $ \Gamma_\psi = f (\gamma_\psi ) . $
If
$$
I_\Delta:=\int_ {\Gamma_\psi} \frac { | \Delta ^{''} (w) | }{|\Delta^{ '} (w)^ 2|} \cdot  |dw| <\infty $$
then
$$
I_\theta:=\int_ {\gamma_\psi} \frac { | \theta ^{''} (z) | }{|\theta^{ '} (z)^ 2|} \cdot  |dz|
<\infty. $$
\end{theor}
\proof
We put $  w = f(z) = \Delta^{-1} (\theta (z) ) , $  $ z \in  \bC ^+ , $ hence

$$  \theta (z) =   \Delta (f(z) ) . $$
Taking first and second derivatives
gives
$$  \theta^{'} (z) =   \Delta^{'} (f(z) ) f^{'}(z) , 
\quad \text{and} \quad 
\theta^{''} (z) =   \Delta^{''} (f(z) ) f^{'}(z) +  \Delta^{'} (f(z) ) f^{''}(z) $$
Next $  w = f(z) $ gives $dw / dz  =  f^{'}(z),  $ and $|dw| =| dz| |f^{'}(z)| . $ Combining we find
$$ \int_ {\gamma_\psi} \frac { | \theta ^{''} (z) | }{|\theta^{ '} (z)^ 2|} \cdot  |dz| = 
\int_ {\Gamma_\psi} \frac { |  \Delta^{''} (w ) f^{'}(z) +  \Delta^{'} (w ) f^{''}(z) | } {|\Delta^{ '} (w)|^2 | f^{'}(z)|^2  } \cdot  \frac{|dw|}{ |f^{'}(z)|}  $$
Applying triangle inequality, we find the following upper bound for the right hand side integral
$$\int_ {\Gamma_\psi} \frac { |  \Delta^{''} (w )|  } {|\Delta^{ '} (w)|^2  } \cdot  \frac{|dw|}{ |f^{'}(z)|}  + \int_ {\Gamma_\psi} \frac {  |f^{''}(z) | } {|\Delta^{ '} (w)| | f^{'}(z)|^3  } \cdot  {|dw|} =
I_1+ I_2  $$

Next recall that $\gamma_\psi $ is a compact subset of $\bC \setminus (-\infty , -\lambda ) , $  and that $f $
is analytic (conformal) with $f^{'} (z ) \ne 0 $ for $z \in \gamma_\psi  $, see Lemma \ref{l18mar-1}.
Hence
\begin{equation}\label{eq17mar-1}
\sup \left\{ \frac{1}{  |f^{'} (z ) | } : z \in \gamma_\psi \right\} \le  C_2,\quad \sup \left\{ \frac{|f''(z)|}{  |f^{'} (z ) |^3 } : z \in \gamma_\psi \right\} \le  C_3. 
\end{equation}
Thus for the first integral we have
$I_1\le C_2 I_\Delta$.

To estimate the second integral $I_2$
we note that
$$
\left(\frac{1}{  \Delta^{'} (w )}\right) ^{'} = -\frac{ \Delta^{''} (w )}{  \Delta^{'2} (w )}.
$$
Hence
\begin{align*}
\frac{1}{  |\Delta^{'} (w )|} =& \left|   \int^{\Delta^{-1}(-\psi +ia)}_w \frac{ \Delta^{''} (v )}{  \Delta^{'2} (v )} dv  +    \frac{1}{  \Delta^{'} (\Delta^{-1}(-\psi +ia) )} \right|   
\\
\le  &
\int^{\Delta^{-1}(-\psi +ia)}_w \frac{ |\Delta^{''} (v )|}{ | \Delta^{'2} (v )|} |dv| 
+   \left | \frac{1}{  \Delta^{'} (\Delta^{-1}(-\psi +ia) )} \right| 
\\
\le  &
\int_{\Gamma_\psi}\frac{ |\Delta^{''} (v )|}{ | \Delta^{'2} (v )|} |dv| 
+   \left | \frac{1}{  \Delta^{'} (\Delta^{-1}(-\psi +ia) )} \right| =I_\Delta+C_4.
\end{align*}
Therefore, see \eqref{eq17mar-1} and  \eqref{eq17mar-3},
\begin{align*}
 I_2 \le&
 C_3 \int_ {\Gamma_\psi} \frac { 1} { |\Delta^{ '} (w)|   } \cdot  {|dw|}=
C_3 \int_ {s_\psi =\Delta^{-1}(\Gamma_\psi)} \frac { 1} { |\Delta^{ '} (w)|^2   } \cdot  {|d\Delta|}
\\
\le& C_3 (I_\Delta+C_4)^2 \int_ {s_\psi}  {|d\Delta|}= C_3 (I_\Delta+C_4)^2 a.
\end{align*}
\endproof

\subsection{Third Reduction}
Let $ n \in \bN .$   Let $ \Gamma_\psi = f (\gamma_\psi ) . $  Define  
$$  \Gamma_{\psi, n} = \{w \in \bC_+ :  \Delta(w) = - \psi +ih : h \in (a/2^{n +1}, a/2^{n }) \} .$$
Let $ \rho^ n \cdot \Gamma_{\psi, n} = \{\rho^ n w \in \bC_+ :  w \in  \Gamma_{\psi, n} \} .$
Recall that $ \Delta (  \rho w ) = 2  \Delta ( w ) , $ for $ w \in  \bC_+  $ and hence
$$ \Delta (  \rho ^n w ) = 2^n  \Delta ( w ) , \quad  w \in  \bC_+  ,$$
which yields the following identity
$$ \rho^ n \cdot \Gamma_{\psi, n} = \{w \in \bC_+ :  \Delta(w)
= - 2^n \psi +ih : h \in (a/2, a) \} .$$

  \begin{theor} \label{t19-3}
$$
\int_ {\Gamma_{\psi ,n}} \frac { | \Delta ^{''} (w) | }{|\Delta^{ '} (w)^ 2|} \cdot  |dw| =
  \frac{ 2^n}{\rho ^n} \int_ { \rho ^n \cdot \Gamma_{\psi ,n} }\frac { | \Delta ^{''} (\omega) | }{|\Delta^{ '} (\omega)^ 2|} \cdot  |d\omega |
 $$
\end{theor}

\proof
Put first $ \omega = \rho ^ n w . $ Then 
$$ \Delta (\omega) = \Delta (\rho^n  w ) = 2^n \Delta (w) . $$
Taking derivatives with respect to $w$ we find
$$\frac{d}{dw} (  \Delta (\rho^n  w ) ) = \Delta^{'} ( \rho^n w )  \rho^n , $$
and
$$
\frac{d}{dw} (2^n \Delta (w) )  =  2^n \Delta^{'} (w) , $$
Since the left hand sides of the above equations coincide, we find that  $\Delta^{'} ( \rho^n w )  \rho^n  =  2^n \Delta^{'} (w)  $  or
$$  \Delta^{'} ( \omega )  = \Delta^{'} ( \rho^n w ) =  \frac{2^n}{  \rho^n } \Delta^{'} (w) ,  $$
Differentiating the above identity again with respect to $ w $ gives
$$ \rho^n  \Delta^{''} ( \rho^n w ) =   \frac{d}{dw} ( \Delta^{'} ( \rho^n w ) ) =   \frac{2^n}{  \rho^n }    \frac{d}{dw}  (\Delta^{'} (w))  =  \frac{2^n}{  \rho^n }  \Delta^{''} ( w ), $$
or 
$$ \Delta^{''} ( \omega) =  \frac{2^n}{  \rho^{2n }}  \Delta^{''} ( w ) .  $$
Since moreover the change of variables  $ \omega = \rho ^ n w  $  yields $|d\omega | = \rho ^n |dw | , $ we finally obtain
$$ \int_ { \rho ^n \cdot \Gamma_{\psi ,n} }\frac { | \Delta ^{''} (\omega) | }{|\Delta^{ '} (\omega)^ 2|} \cdot  |d\omega |
= \int_ {\Gamma_{\psi ,n}} \frac {  \frac{2^n}{  \rho^{2n} }  | \Delta ^{''} (w) | }{  \frac{2^{2n}}{  \rho^{2n} }   |\Delta^{ '} (w)^ 2|} \cdot {  \rho^n }  |dw| =
\frac{\rho^{n}}{{2^n}} \int_ {\Gamma_{\psi ,n}} \frac{ | \Delta ^{''} (w) | }{     |\Delta^{ '} (w)^ 2|} \cdot  |dw|.$$
\endproof

\subsection{Change of Variables}

  Let $ n, N  \in \bN .$  Fix $\psi \in F_N $ with $ 1/2 < \psi < 1 . $
  Recall that we put  
$$  \Gamma_{\psi, n} = \{w \in \bC_+ :  \Delta(w) = - \psi +ih : h \in (a/2^{n +1}, a/2^{n }) \} , $$
 and that in view of
 $$ \Delta (  \rho ^n w ) = 2^n  \Delta ( w ) , \quad  w \in  \bC_+  ,$$
we have 
$$ \rho^ n \cdot \Gamma_{\psi, n} = \{w \in \bC_+ :  \Delta(w)
= - 2^n \psi +ih : h \in (a/2, a) \} .$$
Let $k_{n}$ denote the integer part of  $ 2^n \psi . $ If $k_n $ is odd  we  define $m_n $ by $2m_n = k_n -1$, if $k_n $ is even we  set  $2m_n = k_n $.
A moments reflection shows that
$$
 \rho^ n \cdot \Gamma_{\psi, n} \sbe \cD_{m_n . } 
 $$
 
 Recall that in the definition of $F_N$ we used the dyadic expansion   $\psi = \sum_{k = 1}^\infty \varepsilon_k 2^{-k}$. Therefore the fractional part of $2\psi$ is of the form
 $$
 \underline{s}\psi = \sum_{k = 1}^\infty \varepsilon_{k+1} 2^{-k}.
 $$
 In these notations
 $$
 2^n\psi=k_n+\underline{s}^n\psi.
 $$
 That is, the multiplication acts as the standard shift on the representing sequence $\varepsilon=(\varepsilon_k)$.
 Evidently, this \textit{operation preserves $F_N$}! This observation is highly important in what follows.
 
 \begin{theor} \label{t19-2} There exists $ C < \infty $ such that for each $n \in \bN $
$$ \int_ { \rho ^n \cdot \Gamma_{\psi ,n} }
\frac { | \Delta ^{''} (\omega) | }{|\Delta^{ '} (\omega)^ 2|}\cdot  |d\omega |
\le \frac{C}{ |F^\prime (a_{2m_n}) | }. 
 $$
\end{theor}

\proof
For a change of variables  fix $f_{m_n} : \bC _+ \to \cD_{m_n} , $ and set $ L_n = f_{m_n}^{-1}  (\rho ^n \cdot \Gamma_{\psi ,n} ) $.
Recall, 
$\hat \theta (z) = \theta (P(z)) , z \in \bC_+ $. Thus
\begin{equation}\label{eq21-1}
\Delta (f_{m_n} (z) ) = \hat \theta (z) -2m_n , \qquad z \in \bC_+. 
\end{equation}
%
Since  $2^n\psi=k_n+\underline{s}^n \psi$, we have $L_n =L_{\underline{s}^n\psi}$ if $k_n=2m_n$  and $L_n =L_{\underline{s}^n\psi+1}$ if $k_n=2m_n+1$  as it was defined in Lemma \ref{l21-1}. Due to the invariance property of the shift $\underline{s}$, 
$$L_n\subset K=K(N)$$ for all $n \in \bN$, see the second claim of this lemma. 

Differentiating \eqref{eq21-1} twice we obtain

$$\Delta^\prime (f_{m_n} (z))f^\prime_{m_n} (z) = \hat \theta^\prime  (z) , $$
and
$$ \Delta^{\prime \prime} (f_{m_n} (z))f^{\prime 2}_{m_n} (z) +\Delta^\prime (f_{m_n} (z))f^{\prime \prime} _{m_n} (z)   = \hat \theta^{\prime \prime}  (z) , $$
or
$$ \Delta^{\prime \prime} (f_{m_n} (z))   = \left[\hat \theta^{\prime \prime}  (z)  - \hat \theta^\prime  (z) \frac{f^{\prime \prime} _{m_n} (z)}{f^{\prime }_{m_n} (z)} \right] \frac{1}{f^{\prime 2}_{m_n} (z)} , $$
Consequently with $ \omega =  f_{m_n} (z) , $  and $ L_n = f_{m_n}^{-1}  (\rho ^n \cdot \Gamma_{\psi ,n} ) $ we find

$$ \int_ { \rho ^n \cdot \Gamma_{\psi ,n} }
\frac { | \Delta ^{''} (\omega) | }{|\Delta^{ '} (\omega)^ 2|}\cdot  |d\omega | = \int_{L_n} \frac { | \Delta ^{''} (f_{m_n} (z)) | }{|\Delta^{ '} (f_{m_n} (z))^ 2|} |f^\prime _{m_n} (z) | \cdot |dz| =
\int_{L_n} \frac { \left |\hat \theta^{\prime \prime}  (z)  - \hat \theta^\prime  (z) \frac{f^{\prime \prime} _{m_n} (z)}{f^{\prime }_{m_n} (z)} \right| }{|\hat \theta^{ '}  (z)^ 2|} |f^\prime _{m_n} (z) | \cdot |dz| . $$

Cearly triangle inequality yields that the integral on the right hand side is bounded by

\begin{equation*}
\int_{L_n} \frac { \left |\hat \theta^{\prime \prime}  (z)   \right| }{|\hat \theta^{ '}  (z)^ 2|} |f^\prime _{m_n} (z) | \cdot |dz|  +
\int_{L_n} \frac { \left | \hat \theta^\prime  (z) \frac{f^{\prime \prime} _{m_n} (z)}{f^{\prime }_{m_n} (z)} \right| }{|\hat \theta^{ '}  (z)^ 2|} |f^\prime _{m_n} (z) | \cdot |dz| .
\end{equation*}
By canceling terms in the right hand side integral yields finally the upper bounds
\begin{equation*}
\int_ { \rho ^n \cdot \Gamma_{\psi ,n} }
\frac { | \Delta ^{''} (\omega) | }{|\Delta^{ '} (\omega)^ 2|}\cdot  |d\omega |  \le
\int_{L_n} \frac { \left |\hat \theta^{\prime \prime}  (z)   \right| }{|\hat \theta^{ '}  (z)^ 2|} |f^\prime _{m_n} (z) | \cdot |dz|  +
\int_{L_n} \frac { \left |{f^{\prime \prime} _{m_n} (z)}  \right| }{|\hat \theta^{ '}  (z)|}                       \cdot |dz| .
\end{equation*}

By compactness there exist $ C_1 , C_2 < \infty $ such that 
$$   \sup \{ |\hat \theta^{ '}  (z)|^{-1} :z  \in K \}  \le C_1 $$
and
$$   \sup \{ |\hat \theta^{ \prime \prime }  (z)| |\hat \theta^{ '}  (z)|^{-1} : z \in K \}  \le C_2 .$$
Consequently 
\begin{equation*}
\int_ { \rho ^n \cdot \Gamma_{\psi ,n} }
\frac { | \Delta ^{''} (\omega) | }{|\Delta^{ '} (\omega)^ 2|}\cdot  |d\omega |  \le
C _1 \int_{L_n}   |f^\prime _{m_n} (z) | \cdot |dz|  +
C_2 \int_{L_n} \left |{f^{\prime \prime} _{m_n} (z)}  \right |                       \cdot |dz| .
\end{equation*}
In the next step of the proof we exploit an integral representation for $ f = f_{m_n} $ to show that the right hand side integrals are bounded by $ C |f^\prime (\xi )| ,  $
where $C< \infty$ is independent of $ n \in \bN $, (but it  depends on the compact set $K=K(N)$).

Since 
\begin{enumerate}
\item 
    $f_{m_n} :  \bC_+ \to \bC_+ $ is conformal  ($\Im f_{m_n} (z) > 0 $ for $ z\in \bC_+ $) 
  \item  $f_{m_n}(s) \in \bR$  for $s \in [-\lambda , \lambda]$ , ($\Im f_{m_n} (s) = 0 $
    for $s \in [-\lambda , \lambda]$)
\end{enumerate}    
there exists a non-negative measure $\nu $ supported on $ \bR \setminus [-\lambda , \lambda]$ satisfying
$\int_\bR \frac{d\nu(s)}{1 + s^2} < \infty $  such that
$$ f^\prime_{m_n} (z) = \int_{\bR \setminus [-\lambda , \lambda]}  \frac{1}{(s-z)^2} d\nu(s) $$
  and $$ f^{\prime \prime}_{m_n} (z) = 2\int_{\bR \setminus [-\lambda , \lambda]}  \frac{1}{(s-z)^3} d\nu(s) . $$
    Hence by compactness of $ K \sb \bC_+$
    there exist constants $C_4 = C_4(K) $ such that
    $$
    \sup \{  |f^\prime_{m_n} (z)|  ,   |f^{\prime \prime}_{m_n} (z)| : z \in K \} \le  C_4  |f^\prime _{m_n}(\xi )| . $$
In summary since $ L_n \sbe K , $ there exist $ C_5 = C_5 ( K )$ and  $ C_6 = C_6 ( K )$ such that 
    \begin{equation*} 
    \int_{L_n}   |f^\prime _{m_n} (z) | \cdot |dz|  +
    \int_{L_n} \left |{f^{\prime \prime} _{m_n} (z)}  \right |\cdot |dz|  \le
    C_5\int_{L_n}|dz|  |f^\prime _{m_n}(\xi )|  \le C_6 |f^\prime_{m_n} (\xi )| .
\end{equation*}

Finally, $F(f_{m_n}(z))=P(z)$ and by definition \eqref{eq17mar-5} $f_{m_n}(\xi)=a_{2m_n}$. Therefore
$$ 
F^\prime (f_{m_n} (\xi))f^\prime _{m_n} (\xi)  = P^\prime (\xi)  = 2 \xi . 
$$
That is,
$$
f^\prime _{m_n} (\xi) =\frac{2\xi}{F'(a_{2m_n})}.
$$

\endproof
\subsection{Estimating $|F^\prime(a_{2m_n})|$ from below.}
%

We  review the behaviour of $ f_1 : \bC_+ \to \cD_1 $

\begin{lemma}\label {l20-1}
\begin{enumerate}
\item $z \in \bR$ and $\Im f_1(z) > 0 $ is equivalent to $z \in ( -\infty , -(\lambda^2 -\lambda)) \cup (\lambda , \infty )  $
\item The Herglotz-Nevanlinna measure of $f_1$ is supported in  $ ( -\infty , -(\lambda^2 -\lambda)) \cup (\lambda , \infty )  $
\item $f_1 $ and $ f_1^{\prime} $ restricted to $ [-\xi , \xi ]$  are  continuous.
\end{enumerate}

\end{lemma}
\proof
 Note  that by definition $ f_1(-\xi ) = b_4$, $f_1(\xi) = a_2 $ and by symmetry $ f_1(0) = c_3$.  Further,
$ f_1(-(\lambda^2 -\lambda)) = c_4$.
Indeed, since   $F(f_1(z)) +\xi = P(z) $ and $ F(c_4)  +\xi = P^{\circ 3 } (0 ) = P( \pm P^{\circ 2 }(0) ) = P(  \pm(\lambda^2 -\lambda)) $
we get $ f_1(-(\lambda^2 -\lambda)) = c_4 . $

\begin{lemma} \label{l18mar-5} $$ F(\rho^n f_1(z) ) + \xi =  P^{\circ (n+1) } (z )$$ \end{lemma}
\proof
Set $ \varphi (z) = z + \xi  $ and ${\hat P } (z) = \varphi^{-1} P \varphi (z) . $
Then $\hat P(0) = 0 $ and  ${\hat P } (z) = z(z + 2 \xi) . $
Recall first that Poincare's equation asserts that
$$ F(\rho ^n w ) = \hat P ^{\circ n } (F(w)) . $$
Next
$$  \hat P ^{\circ n } (F(w)) = (\varphi ^{-1} P^{\circ n} \varphi ) (F(w)) =P^{\circ n} (F(w)  +\xi) -\xi $$
Combining with Poincare's equation gives
$$  F(\rho ^n w )  + \xi =  P^{\circ n} (F(w)  +\xi) $$
Inserting $ w = f_1(z) $ and invoking that $F(f_1(z))  +\xi= P(z) $ we find
$$   F(\rho ^n f_1(z) )  + \xi =  P^{\circ n} (F(f_1(z))  +\xi) =   P^{\circ n} (P(z))
=P^{\circ (n+1)} (z)  .   $$

\begin{lemma}
  $$ c_{2^n} = \rho^n c_1, \quad  b_{2^n} = \rho^n b_1,\quad  a_{2^n} = \rho^n a_1 . $$
\end{lemma}
\proof
The proof proceeds by comparing the cuts  of the domains $\Pi _F $ and $\Pi $ at the base points $ - 2^ n . $
Recall that we let $ \Delta $ denote the conformal map from $\bC_+$ onto   $\Pi $.
Let $ (\tilde{a}_{2^ n} ,  \tilde{b}_{2^ n} ) \sb \bR $ denote the pre-image of the cut at  base points $ - 2^ n $ under  the map   $ \Delta $.
Thus
$$  (\tilde{a}_{2^ n} ,  \tilde{b}_{2^ n} ) =  \Delta ^{-1} (\text{cut of  $\Pi $ at  the base point $ - 2^ n $} ) .$$
Next set 
$$  \tilde{c}_{2^ n}  =  \Delta ^{-1} (\text{tip of the cut of  $\Pi $ at  the base point $ - 2^ n $} ) .$$
We have
$$\Delta  (\tilde{c}_{2^ n} ) = - 2^ n +i2^ na =  2^ n( -1 +ia ) =  2^ n\Delta  (\tilde{c}_{1} ) = \Delta  (\rho ^n \tilde{c}_{1} ). $$
Consequently $\tilde{c}_{2^ n}  = \rho ^n \tilde{c}_{1} . $

Similarily 
$$\Delta  (\tilde{a}_{2^ n} ) = - 2^ n  =   2^ n\Delta  (\tilde{a}_{1} ) = \Delta  (\rho ^n \tilde{a}_{1} ) $$
yielding $\tilde{a}_{2^ n}  = \rho ^n \tilde{a}_{1} . $ The same calculation gives  $\tilde{b}_{2^ n}  = \rho ^n \tilde{b}_{1} . $

Recall that we let $ \Psi_F $ denote the conformal map from $\bC_+$ onto   $\Pi_F $.
$$  ({a}_{2^ n} ,  {b}_{2^ n} ) =  \Psi_F ^{-1} (\text{cut of  $\Pi_F $ at  the base point $ - 2^ n $} ) .$$
Next set 
$$  {c}_{2^ n}  =  \Psi_F  ^{-1} (\text{tip of the cut of  $\Pi_F $ at  the base point $ - 2^ n $} ) .$$
We will show next that
$$\tilde{c}_{2^ n} = {c}_{2^ n},  \quad \tilde{a}_{2^ n} = {a}_{2^ n},  \quad \tilde{b}_{2^ n} = {b}_{2^ n} . $$
To this end recall that $ a_{2m}, b_{2m +2} \in \cD_m $ are  zeros of $F ,  $ and that
for $ w \in \cD_m  $
$$ \Delta (w) = \theta (F(w) + \xi )  -2m. $$
Since $F(a_{2^n}) = 0 $ and $\theta (\xi ) = 0 $  we obtain
$$ \Delta (a_{2^n}) =  \theta (F(a_{2^n}) + \xi ) -2^n = \theta (0 +\xi )  -2^n =  -2^n. $$
Since, we defined $\tilde a_{2^n}$ to satisfy  $\Delta (\tilde a_{2^n}) =   -2^n , $ it follows that
$$\Delta (\tilde a_{2^n}) =  \Delta ( a_{2^n}) , $$
and hence $\tilde a_{2^n} =   a_{2^n} .$ We remark that 
$a_{2^n} \in \cD_{2^{n-1}}  . $  Analogously we obtain  $\tilde b_{2^n} =   b_{2^n} .$

Finally we turn to proving  $\tilde c_{2^n} =   c_{2^n} .$ To this end
recall that for odd $k$
\begin{equation}
  F(c_{2^m k}) + \xi = \begin{cases}
                         P(0)   \quad \text{if}\ m = 0 ; \\
                         P^{\circ (m+1)} (0)     \quad \text{if}\ m \neq  0 ; \\
                       \end{cases}
\end{equation}
and that  $c_{2^n} \in \cD_{2^{n-1}}   .$
Hence, invoking B\"ottcher's identity, we obtain
$$ \Delta (c_{2^n}) =  \theta (F(c_{2^n}) + \xi ) -2^n = \theta ( P^{\circ (n+1)} (0))  -2^n = i 2^n a    -2^n. $$
On the other hand we defined
$\tilde c_{2^n}$ to satisfy  $\Delta (\tilde c_{2^n}) =   -2^n  +  i 2^n a    .$  Hence it follows that
$$\Delta (\tilde c_{2^n}) =  \Delta ( c_{2^n}) , $$
and  $\tilde c_{2^n} =   c_{2^n} .$
\endproof


\paragraph{Representing $a_{2 m_n} $ in terms of $f_1$.}  Let $ n, N  \in \bN .$  Fix $\psi_* \in F_N $ with $ 1/2 < \psi_* < 1 . $
  Recall that we put  
$$  \Gamma_{\psi, n} = \{w \in \bC_+ :  \Delta(w) = - \psi +ih : h \in (a/2^{n +1}, a/2^{n }) \} , $$
 and that in view of
 $$ \Delta (  \rho ^n w ) = 2^n  \Delta ( w ) , \quad  w \in  \bC_+  ,$$
we have 
$$ \rho^ n \cdot \Gamma_{\psi_*, n} = \{w \in \bC_+ :  \Delta(w)
= - 2^n \psi_* +ih : h \in (a/2, a) \} .$$
Let $k_{n}$ denote the integer part of  $ 2^n \psi_* . $ As before, if $k_n $ is odd  we  define $m_n $ by $2m_n = k_n -1$, if $k_n $ is even  we set  $2m_n = k_n $.
A moments reflection shows that
$$
 \rho^ n \cdot \Gamma_{\psi_*, n} \sbe \cD_{m_n . } 
$$
   
\begin{lemma}  \label{l18mar-3}There exists $ y \in [-(\lambda^2 -\lambda), \lambda]$  such that  $$\rho^{n-2} f_1(y)   = a_{2m_n} .$$
  If  $ y \in [-(\lambda^2 -\lambda), \lambda]$ then  $f_1(y) \in  [c_4 , c_2] $, 
and if  $d \in [c_4 , c_2] $ then  $\rho^{n-2} d  \in [c_{2^{n+1}}, c_{2^ n}] $.
\end{lemma}
\proof
Note that    $c_2 = \rho c_1$  and $c_4 = \rho ^2 c_1 $ hence 
$$ 
\rho^{n-2} c_4 = c_{2^n} 
\le a_{2m_n}  \le  c_{2^{n-1}} = \rho^{n-2} c_2. $$
Separately $f_1 $ maps the interval $ [-\lambda^2 +\lambda, \lambda] $ onto the interval
$[c_4 , c_2 ] $, hence
$$
{f_1}_{|  [-\lambda^2 +\lambda, \lambda]} = [c_4 , c_2 ] \quad\text{and}\quad 
\rho^{n-2}{f_1}_{|  [-\lambda^2 +\lambda, \lambda]} = [\rho^{n-2}c_4 , \rho^{n-2}c_2 ] .
$$
Since $  a_{2m_n} \in  [\rho^{n-2}c_4 , \rho^{n-2}c_2 ] , $ we verified that there exists
$ y \in  [-\lambda^2 +\lambda, \lambda] $ such that
$$\rho^{n-2}{f_1} (y) = a_{2m_n}  . $$
\endproof

Choose $ \eta > 0 $ such that $P(\eta ) = -\xi . $ If $ \xi > 1 + \sqrt{2} , $
  then $\eta > 1 . $  
Given $n\in \bN  ,$ we  denote by  $P^{\circ n} $ the $n-$fold iterate of $P$. It is a polynomial of degree $ 2 ^ n . $  We let $P^{\circ n \prime} $ denote the complex derivative of $P^{\circ n} .$
\begin{lemma} \label{l18mar-7} Let $n\in \bN  $ and $ y \in  \bC . $  If $P^{\circ n} (y) = \xi , $  then     $ |P^{\circ n \prime} (y) | \ge 2^ n \eta^ n $
\end{lemma}
\proof Since  $P^{\prime } (z)  = 2 z  $  and   $P^{\circ n } (z) = P(P^{\circ (n-1) } (z)) $, applying the chain rule gives,
$$P^{\circ n  \prime } (z) = P^\prime (P^{\circ (n-1) } (z) )P^{\circ (n-1)  \prime } (z) =  2P^{\circ (n-1) } (z)P^{\circ (n-1)  \prime } (z) . $$
Using the above identity as  the basis of the recursion, we obtain the following product representation for $P^{\circ n  \prime } (z)$, 
$$P^{\circ n  \prime } (z) = 2^n P^{\circ (n-1)  } (z) P^{\circ (n-2)  } (z)\cdots  P (z) z . $$
Next,  since   $ P^{\circ k}(P^{\circ (n-k)} (y)) = \xi  $,  we find that 
$$P^{\circ (n-k)} (y)  \in E_0\sb [- \xi,  \xi] \setminus (-\eta, \eta) ,$$
for each $1 \le k \le n-1$ and hence that
$|P^{\circ (n-k)} (y)| \ge \eta . $ Inserting this into the above product representation for $P^{\circ n  \prime } (z)$ we obtain
$$ |P^{\circ n\prime} (y) |  \ge 2^ n \eta^ n . $$
\endproof

  \begin{lemma}\label{l19-1} There exists $C$ such that for all $n\in\bN$
  \begin{equation}\label{eq19mar-7}
  |F'(a_{2m_n} )|^{-1} \le C \frac {\rho^n }{2^n\eta^n } . 
    \end{equation}
  \end{lemma}
  \proof
  By Lemma \ref{l18mar-3}, and Lemma \ref{l18mar-5} 
  $$ F(a_{2m_n} ) + \xi =   F(\rho^{n-2} f_1(y) ) + \xi =  P^{\circ (n-1) } (y )$$
  Since $ F(a_{2m_n} ) = 0 $ we infer that $ P^{\circ (n-1) } (y ) = \xi . $ By Lemma \ref{l18mar-7}
  \begin{equation}\label{eq18mar-7}
  |(P^{\circ (n-1) })^\prime (y )| \ge  2^{n-1} \eta^{n-1} . 
  \end{equation}

  Next invoke the identity  $F(\rho^{n-2} f_1(z) ) + \xi =  P^{\circ (n-1) } (z ) $ of Lemma \ref{l18mar-5}
  and differentiate both sides with respect to $z$. This yields
  $$ F^\prime(\rho^{n-2} f_1(z) ) \rho^{n-2} f_1^\prime (z) =   P^{\circ (n-1) \prime} (z ) , $$
  or
  $$ \frac{1}{F^\prime(\rho^{n-2} f_1(z) )} = \frac{ \rho^{n-2} f_1^\prime (z)}{   P^{\circ (n-1) \prime} (z )} , $$
  Specializing to $ z = y$ and invoking  that $a_{2m_n}  =   \rho^{n-2} f_1(y)$ and \eqref{eq18mar-7}   gives
  $$ \frac{1}{|F^\prime(a_{2m_n} )|} \le  \frac{\rho^{n-2} |f_1^\prime (y)|}{2^{n-2} \eta^{n-2}} , $$
  Since $y \in  [-\xi , \xi ] $ and by Lemma \ref {l20-1} $\sup_{ x \in [-\xi , \xi ] } | f^\prime _{1} (x)|$ is finite, we get  \eqref{eq19mar-7}.
  \endproof

  \subsection{Proof of the Main Theorem \ref{main}}
  Let $P(z)$ be an even monic quadratic polynomial with a fixed point  $\xi>1+\sqrt{2}$, $P(\xi)=\xi$. Equivalently, $\eta>1$, where $P(\eta)=-\xi$.
  
  We assume that there exists $N\in\bN$ such that $\psi\in F_N$. WLOG $\psi\in (1/2,1)$ (otherwise we consider $1-\psi\in F_N$, see below). For an arbitrary $n\in\bN$ we have
  $$
  2^n\psi=k_n+\underline{s}^n\psi, \quad \underline{s}^n\psi\in F_N,
  $$
  and we let $k_n=2m_n$ if $k_n$ is even and $k_n-1=2m_n$ if $k_n$ is odd.
  
  Due to Lemma \ref{l19-1} uniformly in $n\in\bN$ there exists $C_1$ such that the derivative of the Poincare function $F$ can be estimated from below at its zero $a_{2m_n}$:
  $$
  \frac 1 {|F'(a_{2m_n})|}\le C_1\frac{\rho^n}{2^n\eta^n}
  $$
  
  By Theorem \ref{t19-2}, there exists $C_2$ such that the integral on the path $\rho^n\cdot \Gamma_{\psi,n}$ of the conformal mapping $\Delta:\bC_+\to \Pi$ can be estimated  as
  $$
\int_ { \rho ^n \cdot \Gamma_{\psi ,n} }
\frac { | \Delta ^{''} (\omega) | }{|\Delta^{ '} (\omega)^ 2|}\cdot  |d\omega |\le C_2 \frac{\rho^n}{2^n\eta^n}.
  $$
Respectively, by Theorem \ref{t19-3}
$$
\int_ {\Gamma_{\psi ,n}} \frac { | \Delta ^{''} (w) | }{|\Delta^{ '} (w)^ 2|} \cdot  |dw| \le C_2\frac{1}{\eta^n}
$$
  and therefore
  $$
I_\Delta=\int_ {\Gamma_{\psi }} \frac { | \Delta ^{''} (w) | }{|\Delta^{ '} (w)^ 2|} \cdot  |dw| =
\sum_{n=0}^\infty\int_ {\Gamma_{\psi,n }} \frac { | \Delta ^{''} (w) | }{|\Delta^{ '} (w)^ 2|} \cdot  |dw| 
\le C_2\frac{\eta}{\eta-1}<\infty.
$$

By Theorem \ref{t19-4}
$$
I_\theta:=\int_ {\gamma_\psi} \frac { | \theta ^{''} (z) | }{|\theta^{ '} (z)^ 2|} \cdot  |dz|
<\infty$$
and by Theorem \ref{t19-5}
$$ 
\int _{ R_{\psi}} |T^{''} (w) | \cdot |d w | <\infty.
$$
Due to the symmetry of the Green's domain $H$
$$ 
\int _{ R_{1-\psi}} |T^{''} (w) | \cdot |d w | =\int _{ R_{\psi}} |T^{''} (-\bar w) | \cdot |d (-\bar w) |<\infty.
$$
Therefore the initial restriction $\psi\in(1/2,1)$ can be removed.

The second statement of Theorem \ref{main} was proved in Proposition \ref{p19-6}.

\bibliographystyle{abbrv} 
\bibliography{hardymartingales}
{\bf  Authors Addresses:}\\ 
Department of Mathematics\\
University of Toronto\\
Toronto  M5S 2E4 \\
ilia@math.toronto.edu \\ \\
Department of Mathematics\\
J. Kepler Universit\"at Linz\\
A-4040 Linz\\
paul.mueller@jku.at\\ \\
 Department of Mathematics\\
J. Kepler Universit\"at Linz\\
A-4040 Linz\\
peter.yuditskii@gmail.com
\end{document}